\pdfoutput=1
\documentclass[11pt]{article}
\usepackage[utf8]{inputenc}
\usepackage{mystyle}

\usepackage{graphicx}

\graphicspath{{./figures/}}

\title{A $\Gamma$-Convergence Result for the Upper Bound\\Limit Analysis of Plates}
\author{Jeremy Bleyer \footnotemark[1]
\and Guillaume Carlier \footnotemark[2]
\and Vincent Duval \footnotemark[3]
\and Jean-Marie Mirebeau \footnotemark[4]
\and Gabriel Peyr\'e \footnotemark[4]
}

\usepackage{hyperref}

\begin{document}
\date{}
\maketitle
\footnotetext[1]{Universit\'e Paris-Est, Laboratoire Navier, Ecole des Ponts ParisTech-IFSTTAR-CNRS (UMR 8205) (\url{jeremy.bleyer@enpc.fr})} 
\footnotetext[2]{CEREMADE, Universit\'e Paris-Dauphine (\url{carlier@ceremade.dauphine.fr})} 
\footnotetext[3]{INRIA, Domaine de Voluceau, Rocquencourt (\url{vincent.duval@inria.fr})} 
\footnotetext[4]{CNRS, CEREMADE, Universit\'e Paris-Dauphine (\url{{mirebeau,peyre}@ceremade.dauphine.fr})} 
\renewcommand{\thefootnote}{\arabic{footnote}} 


\begin{abstract}
Upper bound limit analysis allows one to evaluate directly the ultimate load of structures without 
performing a cumbersome incremental analysis. In order to numerically apply this method to thin plates in bending, several authors have proposed to use various finite elements discretizations. We provide in this paper a mathematical analysis which ensures the convergence of the finite element method, even with finite elements with discontinuous derivatives such as the quadratic 6 node Lagrange triangles and the cubic Hermite triangles. More precisely, we prove the $\Gamma$-convergence of the discretized problems towards the continuous limit analysis problem. Numerical results illustrate the relevance of this analysis for the yield design of both homogeneous and non-homogeneous materials. 
\end{abstract}

{\noindent\textbf{Keywords:~}{Bounded Hessian functions; Finite Element Method; $\Gamma$-convergence.}\\
~\\
\noindent\textbf{AMS Subject Classification:~}{ 74K20, 74S05, 74G15}
}



\section{Introduction and Motivation of the Model}

Yield design or limit analysis theory aims at evaluating directly the ultimate load of mechanical structures using the compatibility between equilibrium equations and a local strength criterion. More precisely, assuming that the structure $\Omega$ is subject to a multiplicative loading $\lambda L$, the ultimate load multiplicative factor $\lambda^+$ is obtained by solving the following formal problem
\[ 
	\lambda^+ = \max_{\lambda\in\mathbb{R},\Sigma \in \Ss} \lambda 
	\qstq
	\forall x\in \Omega, \quad
	\choice{
		\Ee(\Sigma)(x) = \lambda L(x),  \\
		\Sigma(x) \in G(x), 
	}
\]
where $\Sigma(x)$ represents the internal forces at a point $x$ of the structure $\Omega$ and where $\Ss$ collects all fields $\Sigma$ of internal forces which satisfy the stress boundary conditions, $\Ee$ is a linear differential operator corresponding to the equilibrium equations and $G(x)$ is a convex set representing the strength criterion at a point $x$ in the structure $\Omega$. This means that the ultimate load represents the greatest load factor such that there exists at least one field of internal forces in the structure satisfying both the equilibrium equation and the strength criterion at each point. 

The dual problem can be expressed formally as follows :
\begin{align}\label{eq-intro-liman}
	\lambda^+ = \min_{u \in \Cc} \int_{\Omega}\pi(x,\Dd u)dx 
  	\qstq
	\int_{\Omega} L(x)u(x) dx=1 
\end{align}
where $\Dd  = \Ee^*$ is the dual equilibrium operator or the strain operator, $u$ is a virtual velocity field which has to satisfy the velocity boundary conditions and $\pi(x,\cdot)$ is the support function of $G(x)$. This dual formulation means that the ultimate load corresponds to the minimum of the maximal resisting work over all virtual velocity fields satisfying a normalization condition on the work of external forces. Both approaches, the primal one (known as the static approach) and the dual one (known as the kinematical approach), are equivalent and can be implemented numerically by an appropriate discretization of the stress or velocity fields, yielding to a lower bound or an upper bound of the ultimate load factor.

In this work, we are interested in the kinematical approach~\eqref{eq-intro-liman} with virtual velocity fields for a specific mechanical model, namely the thin plate model. In this model, $\Omega$ is a region of $\mathbb{R}^2$ and the internal forces are the bending moments and are represented by a symmetric matrix field $M(x)$ (see~\cite{ciarlet1979plate} for the derivation of this model). The equilibrium operator is a second-order operator namely $\Ee(M) = \textrm{div\,div\,}M$. For this reason, the strain rate (or curvature rate for the plate model) is given by the Hessian matrix of the scalar velocity field $u$ as $\kappa = \Dd u = D^2 u$.

\subsection{Previous Works}

\paragraph{Numerical limit analysis of thin plate in bending problems.} 

From a numerical point of view, the resolution of such limit analysis problems is traditionally performed by discretizing the structure into finite elements while the corresponding discrete variational inequality problem is formulated as a convex programming problem. Early works considered a linear interpolation of the velocity in each triangular finite element \cite{munro1978yield}. The corresponding mathematical problem can be formulated as a linear programming problem and can be solved using simplex or interior-point algorithms. This very simple discretization is known as the yield line method in the field of plate mechanics, since each finite element admits a potential discontinuity of the velocity gradient (corresponding to the plate rotation) along its edges. In particular, no curvature deformation occurs inside any element since the Hessian is identically zero. Therefore, this method fails in general to predict the exact collapse load even with an infinitely refined mesh as pointed out by Braestrup \cite{braestrup1971yield}, the computed upper bound being then very sensitive to the mesh layout \cite{johnson1994mechanism,jennings1996identification}. Hodge and Belytschko considered a quadratic interpolation of the velocity field ensuring only $C^0$ continuity between elements \cite{hodge1968numerical} whereas later works used $C^1$-continuous elements \cite{capsoni1999limit,le2010upper}. More recently, it was suggested in \cite{bleyer2013non-conforming} that the use of finite elements ensuring $C^0$ continuity only exhibits better convergence rates to the exact collapse load than elements ensuring $C^1$-continuity. Finally, one can also mention recent papers relying on a meshfree discretization \cite{le2009limit,zhou2012upper}.

\paragraph{Mathematical models and their analysis.} 

The function $\pi$ appearing in the limit analysis problem~\eqref{eq-intro-liman} having linear growth and the strain  operator $\Dd$ being of second order, the adequate functional space for the study of~\eqref{eq-intro-liman} is the space of \textit{bounded Hessian  functions}, which consists of functions whose second derivative is a bounded Radon measure. 
This space was introduced and studied in the pioneering work of Demengel and Temam~\cite{demengel1983problemes,demengel1984fonctions,temam1985mathematical,demengel1989compactness}: the abstract properties of bounded Hessian functions are established in~\cite{demengel1984fonctions} (see also~\cite{temam1985mathematical}), whereas~\cite{demengel1983problemes} studies the limit analysis problem in the homogeneous case ($\pi(x,\Dd u)={\Pi}(\Dd u)$): existence of solutions and their characterization are stated. Compactness properties are studied in~\cite{demengel1989compactness} with an application to the limit load problem.
Finer mechanical models may also involve the space of \textit{bounded Hessian functions}: in~\cite{hadhri1985-hbbd,hadhri1988}, it is used together with the space of \textit{bounded deformations} to model the bending and the compression of a plate constituted by an elastoplastic material.

The strain $\Dd u$ being a Radon measure, it raises the issue of the definition of $\int_{\Omega}\pi(x,\Dd u)dx$ as a convex function of a measure, which was tackled in \cite{demengel1984convex} in the case of a uniform penalty $\pi$ (independent of $x$) and in \cite{hadhri1985-convexe} in the general case. This extended formulation allows to deal with inhomogeneous penalties as in~\cite{hadhri1986-these}. Eventually Telega  studies in~\cite{Telega1995epilimit} the case of a periodically inhomogeneous plate material, and he shows the $\Gamma$-convergence of the problem towards the homogeneous problem as the period vanishes (see below for the definition of $\Gamma$-convergence).

\subsection{Contributions and Outline of the Paper}

This work is a companion paper to the numerical study~\cite{bleyer2013non-conforming}. We consider a relaxed formulation of~\eqref{eq-intro-liman} (similar relaxations are shown to be tight in~\cite{demengel1983problemes} in the study of the homogeneous case),
\begin{align}
  \inf_{u\in \HB} \int_{\Omega} \pi\left(x,D^2u\right) +\int_{\Gamma_N} \pi\left(x,-\frac{\partial u}{\partial \nu}\nu\otimes \nu\right)d\Hh^1
  \label{eq-intro-pb}
\end{align}
where $\Gamma_N \subseteq \partial \Omega$, and we impose the constraint that $u=0$ on $\Gamma_D\subseteq \partial \Omega$ and that the load $L$ should satisfy $\langle L, u\rangle =1$. The precise definitions of each notion is given in Section~\ref{sec-preliminaries}. As for $\int_{\Omega} \pi\left(x,D^2u\right)$, we  rely on a formulation due to Reshetnyak which is different from~\cite{hadhri1985-convexe} and which is exposed in~\cite{Ambrosio}. Its main advantage is that it provides continuity properties provided $\pi$ is continuous in the first variable. 

Our main contribution is to ensure the consistency of the finite element method used in~\cite{bleyer2013non-conforming} to solve~\eqref{eq-intro-pb}. \textbf{We prove that the discretized problems $\Gamma$-converge towards Problem~\eqref{eq-intro-pb} as the size of the triangulation goes to zero}.
As a consequence every cluster point (which always exists) of the solutions to the discrete problems is a solution to~\eqref{eq-intro-pb}.

To our knowledge, although $\Gamma$-convergence results of finite element approximations exist in the literature for first order strain/equilibrium operators~ (\textit{e.g.}~\cite{bendhia1989fem,Ortner2011fem}), this kind of issue has never been tackled for the second-order operators involved in the thin plate limit analysis. Such a problem raises specific issues like the use of the space of bounded Hessian functions, (\ie whose second derivative is a Radon measure) as opposed to the classical Sobolev spaces $W^{1,p}$ or $W^{2,p}$, $p>1$.

The outline of the paper is as follows. In Section~\ref{sec-preliminaries}, several facts about the space of functions with bounded Hessian are recalled. The hypotheses about the domain and the penalty function $\pi$ are stated precisely.
The formulation of the continuous problem is given in Section~\ref{sec-continuous-pb}. For the sake of completeness, we give a proof of existence of solutions, although the argument is similar to the one used by Demengel~\cite{demengel1983problemes} for the homogeneous case. We also state an approximation result of functions with bounded Hessian with functions which are $C^3$ up to the boundary.
The main result of the paper is stated in Section~\ref{sec-elfini}, where we describe the finite element approximations and we prove their $\Gamma$-convergence towards the continuous problem. Eventually, Section~\ref{sec-numerics} provides numerical experiments which illustrate the efficiency of the approximation including in both the homogeneous (constant $\pi$) and the inhomogeneous cases. 

\paragraph{Notations.}

We adopt the following definitions and notations throughout the paper. 

$\tens$ is the set of real symmetric tensors of order 2. For $A,B\in \tens$, their double inner product is $A\tprod B=\sum_{1\leq i,j\leq 2} A_{ij}B_{ij}$.
The Frobenius norm of $A\in \tens$ defined by $|A|_F = \sqrt{A\tprod A}$ is denoted by $|A|_F$ or simply $|A|$.

We denote by $\Ll^2$ the Lebesgue measure on $\RR^2$, by $\Hh^1$ the one-dimensional Hausdorff measure. Given a Borel set $B\subseteq \RR^2$, $\Hh^1\restr B$ refers to the restriction of $\Hh^1$ to $B$, that is $(\Hh^1\restr B)(E)= \Hh^1(B\cap E)$ for all $E\subseteq \RR^2$. More generally, we write $(\mu\restr B) (E)=\mu(B\cap E)$ for any Borel measure $\mu$ and any $\mu$-measurable set $E\in\RR^2$. Given a finite dimensional vector space $X$, $\Mm_b(\Omega,X)$ will denote the set of finite $X$-valued Radon measures on the open set $\Omega$. For all sets $E\subset \Omega$, $E\subcomp\Omega$ means that its closure $\overline{E}$ is compact and $\overline{E}\subset\Omega$.

Given a locally integrable function $u:\Omega\rightarrow \RR$, we denote by $Du$ (resp. $D^2u$) its distributional derivative (resp. second derivative).
If $Du$ is representable by some locally integrable function (\ie $u$ has a Sobolev regularity), we refer to it as $\nabla u$. A similar convention is adopted for $D^2u$ and $\nabla^2 u$. The restriction of $u$ to some set $B\subseteq \RR^2$ is denoted by $u_{\vert B}$.

Given a real vector space $X$ and a function $q:X\rightarrow [-\infty, +\infty]$, we say that $q$ is positively $1$-homogeneous if for all $x\in X$, and all $t>0$
 $q(tx)=tq(x)$.

Eventually, $\Pp_k$ denotes the space of polynomials of two variables of degree $k$ or less and $\mathbb{S}^{N-1}$ denotes the unit sphere of $\RR^N$.


\section{Preliminaries}
\label{sec-preliminaries}

\subsection{Domain Regularity and Boundary Conditions}
\label{sec-domain-assumption}

Throughout this paper, we shall consider a piecewise $C^3$-regular domain $\Omega\subset \RR^2$. 
More precisely, we assume that $\Omega\subset \RR^2$ is a bounded connected open set and there exists a finite set $\{x_1,\ldots, x_J\}\subseteq \partial \Omega$ such that 
\begin{itemize}
  \item for all $x\in \partial \Omega\setminus \{x_1,\ldots, x_J\}$, $\partial \Omega$ is $C^3$ regular in a neighborhood of $x$,
  \item  there exist pairwise disjoint open neighborhoods $V_j$ of $x_j$, for each $j\in \{1, \ldots, J\}$ and $C^3$-diffeomorphisms $\phi_j$ from $V_j$ to $(-1,1)^2$ such that
  \begin{align*}
    \phi_j(x_j)&=(0,0)\\
  \phi_j(V_j\cap \Omega)&=(0,1)^2.
  \end{align*}
\end{itemize}
In the last section, when studying the convergence of the finite element approximation, we shall make the stronger assumption that $\Omega$ is a convex polytope.

\begin{figure}
  \centering
  \begin{tikzpicture}[scale=2.0]
\definecolor{vert}{rgb}{0.15 , 0.5 , 0.25} 
\definecolor{vert}{rgb}{0.1725 , 0.6627 , 0.1725} 
 \begin{scope}[shift={(0,0)}]
   \coordinate (A) at (-0.3,0.4) ;
  \coordinate (B) at (0.5,-0.7) ;
  \coordinate (C) at (2.1,-0.1) ;
  \coordinate (D) at (1.3,1.1);
  \coordinate (E) at (0.5,1.1);
  \coordinate (O) at (1.0,0.5);
  \coordinate (Ot) at (1.6,0.85);
  \coordinate (Otb) at (0.15,-0.3);

  \coordinate (Abis) at (0.2,0.7);
  \coordinate (Bbis) at (0.9,-0.6);
  \coordinate (Bter) at (1.0,0.4);

  \path[right color=gray!20, right color=gray!80,draw=black!80] (C) parabola bend (1.3,1.3) (E)  -- cycle;
  \path[right color=gray!20, right color=gray!80,draw=black!80] (A) parabola bend (0.1,-0.5) (B)  -- (Abis) -- cycle;
  \path[left color=gray!20, right color=gray!80,draw=black!80] (A) node[left] {$x_1$} to[bend left=10] (B) node[below] {$x_2$} ..controls (Bbis) and  (Bter)..  (C) node[right] {$x_3$} to[bend left=12] (D) node[above] {$x_4$} to[bend right=12] (E) node[above left] {$x_5$} to[bend left=10] (A) -- cycle;

  \path[thick,dashed, draw=blue!80] (A) to[bend left=10] (B);
  \path[thick,dashed, draw=blue!80] (C) to[bend left=12] (D) to[bend right=12] (E) to[bend left=10] (A);
  \coordinate (GD) at (0.25,0.6);
  \draw[color=blue] (GD) node {$\Gamma_D$};  

  \path[thick, draw=blue!80] (A) to[bend left=10] (B);
  \path[thick, draw=blue!80] (C) to[bend left=12] (D) to[bend right=12] (E);
  \coordinate (GN) at (1.2,0.9);
  \draw[color=blue] (GN) node {$\Gamma_N$};  

  \draw[color=black!60] (O) node {$\displaystyle{\Omega}$};
  \draw[color=black!70] (Ot) node {$\displaystyle{\tilde{\Omega}}$};
  \draw[color=black!70] (Otb) node {$\displaystyle{\tilde{\Omega}}$};
  \fill (A)  circle[radius=1pt];
  \fill (B)  circle[radius=1pt];
  \fill (C)  circle[radius=1pt];
  \fill (D)  circle[radius=1pt];
  \fill (E)  circle[radius=1pt];

  \end{scope}

\end{tikzpicture}
  \caption{The domain $\Omega$ is piecewise $C^3$ regular, \ie $C^3$ except at a finite number of points $\{x_1,\ldots x_J\}$. The (relaxed) Neumann boundary condition
  is imposed on $\Gamma_N=(\partial \Omega)\cap \tilde{\Omega}$ where $\tilde{\Omega}\supset \Omega$ is a piecewise $C^3$ domain.
  }
  \label{fig-omega}
\end{figure}
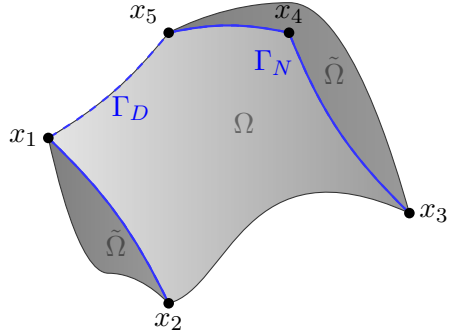

As for the boundary conditions imposed to Problem~\eqref{eq-intro-pb}, we shall assume that the Dirichlet condition $u=0$ (or more precisely $\gamma_0 u=0$, see the definition in Section~\ref{sec-hb-facts}) is imposed on $\Gamma_D$, where $\Gamma_D$ is a union of connected components of $\partial \Omega\setminus \{x_1, \ldots, x_J\}$. Moreover, the Neumann boundary condition  $\Dpart{u}{\nu}=0$ (in fact $\gamma_1 u=0$, see below) is imposed on $\Gamma_N$, where $\Gamma_N$ is a union of connected components of $\partial \Omega\setminus \{x_1, \ldots, x_J\}$ such that $\Gamma_N\subseteq \Gamma_D$. Observe that we work in fact with a relaxed Neumann condition, as shown in the formulation~\eqref{eq-intro-pb}.

From the above assumptions, we deduce that there exists a bounded open connected set $\tilde{\Omega}\subseteq \RR^2$ with piecewise smooth boundary such that $\Omega\subset \tilde{\Omega}$ and $\Gamma_N=(\partial \Omega) \cap \tilde{\Omega}$ (modulo the points $(x_j)_{1\leq j\leq J}$, see Figure~\ref{fig-omega}).

Eventually, we assume that $\Gamma_D$ contains at least three points which are not aligned. 

\subsection{The Penalty Function $\pi$}
\label{sec-pi-hypo}


The mechanical strength of the plate is described by a closed convex set $G\subset \tens$ such that any bending state outside of $G$ is impossible to achieve whereas a bending state lying on the boundary of $G$ means that the plate attained its full strength capacity at this given point. For thin plates in bending, the convex set $G$ is supposed to be bounded and to contain 0. Typical examples of such criteria are the following:
\begin{itemize}
\item the von Mises criterion: it is characterized by a mechanical parameter $M_0>0$ which represents the plate strength in uniaxial bending and it is given by the following quadratic norm on the components $M_{ij}$ of $M$:
\[ G=\{M\in \tens \:|\: \sqrt{M_{11}^2+M_{22}^2-M_{11}M_{22}+3M_{12}^2} \leq M_0\} \] 
\item the Tresca criterion: it is also characterized by a parameter $M_0$ having the same mechanical interpretation but here the criterion is expressed in terms of the eigenvalues $M_I$ and $M_{II}$ of the bending moment $M$ so that it is different from the von Mises criterion: 
\[ G=\{M\in \tens \:|\: \max(|M_I|,|M_{II}|,|M_I-M_{II}|) \leq M_0\} \] 
\item the Johansen criterion: it can be characterized by two parameters $M_0^+,M_0^->0$ which represent the plate strength in positive and negative uniaxial bending and it is given by:
\[ G=\{M\in \tens \:|\: -M_0^-\leq M_I,M_{II} \leq M_0^+\} \] 
\end{itemize} 
 We refer the reader to~\cite{prager1959introduction} for more details about the physical concepts involved in the limit analysis of thin plates.

 It is equivalent to specify $G$ or to specify its support function $\pi$ (see~\cite{rockafellar1986convex}) which appears naturally in the kinematical approach:
\[ \pi(\kappa) = \sup \{M:\kappa\,,\, M\in G\}. \]
It can be shown that $\pi$ is convex continuous and positively homogeneous. The support functions corresponding to the previous classical strength criteria are respectively given by:
\begin{align*}    
\mbox{(von Mises)}\quad  \pi(\kappa) &= \dfrac{2M_0}{\sqrt{3}}\sqrt{\kappa_{11}^2+\kappa_{22}^2+\kappa_{11}\kappa_{22}+\kappa_{12}^2},\\
\mbox{(Tresca)}\quad  \pi(\kappa) &= M_0\max(|\kappa_I|,|\kappa_{II}|,|\kappa_I+\kappa_{II}|),\\
\mbox{(Johansen)}\quad \pi(\kappa) &= \max(M_0^+\kappa_I,-M_0^-\kappa_I) + \max(M_0^+\kappa_{II},-M_0^-\kappa_{II}). 
\end{align*}

Depending on the local width or composition of the plate, the convex set $G=G(x)$ may vary both in size or shape. As a consequence the penalty function may also vary in $\Omega$: $\pi(\kappa)=\pi(x,\kappa)$. 

In the following, we shall consider a continuous function $\pi : \tilde{\Omega}\times \tens\rightarrow [0,\infty)$ which is convex and positively $1$-homogeneous in the second variable. Moreover, 
we assume that there exists $\alpha>0$, $\beta>0$ such that for all $x\in \tilde{\Omega}$ and $\kappa\in \tens$,
\begin{align}
  \alpha |\kappa|_F\leq \pi(x,\kappa)\leq \beta|\kappa|_F.
  \label{eq-bounded-coercive}
\end{align}
Then it can be shown that $\pi(x,\cdot)$ is the support function of some closed convex set $G(x)\subseteq \tens$, and~\eqref{eq-bounded-coercive} means that $G(x)$ is bounded and contains $0$ in its interior (uniformly in $x$).

\subsection{Bounded Hessian Functions}
\label{sec-hb-facts}

We recall here the main properties of the space of functions with bounded Hessian. We refer the reader to \cite{demengel1984fonctions,temam1985mathematical,demengel1989compactness} for more details about this space (also see the textbook ~\cite{Ambrosio}).
A function $u\in W^{1,1}(\Omega)$ whose second derivative $D^2u$ in the sense of distributions is a finite (matrix valued) Radon measure is called a bounded Hessian function. 

The total variation of $D^2u$, denoted by $|D^2u|$, is then a finite positive Radon measure on $\Omega$, and its total mass, given by:
\begin{align}
\sup \bigg\{ \int_\Omega u(x)\div\div \varphi(x)dx;
\   \varphi \in C_c^\infty(\Omega, \tens), \forall x\in \Omega\ |\varphi(x)|_F \leq 1\bigg\},
\label{eq-hb-totvar}\\
\qwhereq \div\div \varphi= \sum_{i,j}\frac{\partial^2\varphi_{i,j}}{\partial x_i\partial x_j},\nonumber
\end{align}
will be denoted by $|D^2u|(\Omega)$ or by $\int_{\Omega} |D^2u|$. One may observe from~\eqref{eq-hb-totvar} that the map $u\mapsto |D^2u|(\Omega)\in [0,+\infty]$ is well-defined on $L^1_{loc}(\Omega)$ and is lower semi-continuous. In fact, $u\in L^{1}(\Omega)$ is a Bounded Hessian function if and only if $|D^2u|(\Omega)<+\infty$, provided that $\Omega$ has the cone property (see~\cite{Adams1975} for the definition, this is obviously the case here since $\Omega$ is piecewise $C^3$).


Endowed with the norm $\normhb{\cdot}$ defined by: 
\begin{align}
  \forall u\in \HB,\quad \normhb{u}\eqdef \|u\|_1 +\|\nabla u\|_1+|D^2u|(\Omega),
\label{eq:norme}
\end{align}
the space $\HB$ of bounded Hessian functions is a Banach space. If $\Omega$ has the cone property, 
$\normhb{\cdot}$ may be replaced with the equivalent norm defined by $\|u\|_1+|D^2u|(\Omega)$ for all $u\in \HB$ (see \cite[Prop. 1.3]{demengel1984fonctions}), and the injection $\HB\hookrightarrow W^{1,1}(\Omega)$ is compact. Moreover, the injection $\HB\hookrightarrow C(\overline{\Omega})$ is continuous.

Two other topologies, known as the ``weak'' and ``intermediate'' topologies, are introduced in~\cite{demengel1984fonctions}. 
The weak topology is the one induced by the norm $\normw{\cdot}$ and the family of seminorms $u\mapsto \left|\int_{\Omega} \psi d(D^2u)_{ij} \right|$ for all $\psi \in C^{\infty}_c(\Omega)$ and $1\leq i,j\leq 2$. Here, $(D^2u)_{ij}$ refers to each entry of $D^2u$, so that it is a signed Radon measure, and in fact $\int_{\Omega} \psi d(D^2u)_{ij}= \int_{\Omega}u(x)\frac{\partial^2\varphi_{i,j}}{\partial x_i\partial x_j}(x)dx$.
An important property is that if a sequence of functions $(u_n)_{n\in\NN}$ is bounded in $\HB$ (for the strong topology), then there exists a function $u_\infty\in \HB$ and a subsequence $(u_{n'})_{n'\in\NN}$ that converges towards $u_\infty$ for the weak topology.

The intermediate topology is defined by the distance
\begin{align}
  d(u_1,u_2)\eqdef \|u_1-u_2\|_{1} + \left| |D^2u_1|(\Omega) -|D^2u_2|(\Omega) \right|.
  \label{eq-hb-distintermed}
\end{align}

Although it is not completely obvious, the intermediate topology is finer than the weak topology, as shown in the next proposition.
\begin{prop}
  Let $u\in \HB$ and $(u_n)_{n\in\NN}\in (\HB)^\NN$ such that 
  $$ \lim_{n\to +\infty}d(u,u_n)=0. $$  
  Then $u_n$ converges towards $u$ for the $\HB$-weak topology.
  \label{prop-intermed-weak}
\end{prop}
\begin{proof}
  The main point to prove is that $\lim_{n\to +\infty} \|\nabla(u-u_n)\|_{L^{1}}=0$.
  
Assume by contradiction that there is some $\varepsilon>0$ and a subsequence (still denoted $u_n$) such that $\|\nabla(u_n-u)\|_{L^1}\geq \varepsilon$. The sequence $\|u_n\|_{L^1}+|D^2u_n|(\Omega)$ being bounded (so that $\|u_n\|_{\HB}$ is bounded as well, by equivalence of the norms) we may extract a subsequence which converges in $W^{1,1}(\Omega)$ towards some $\tilde{u}\in \HB$. In particular we have $\lim_{n\to+\infty}\|\nabla (u_n-\tilde{u})\|_{L^1}=0$. Since on the other hand $\lim_{n\to+\infty} \|u_n-u\|_{L^1}=0$, we must have $\tilde{u}=u$, and thus $\lim_{n\to+\infty}\|\nabla (u_n-u)\|_{L^1}=0$, which contradicts the assumption. 
Thus $\lim_{n\to +\infty} \|\nabla(u-u_n)\|_{L^{1}}=0$.

The weak-* convergence of $D^2u_n$ towards $D^2u$ comes from the fact that $\lim_{n\to +\infty} u_n=u$ in $L^1(\Omega)$ and that $\sup_{n\in\NN}\left(|D^2u_n|(\Omega)\right)<+\infty$.
\end{proof}

The following important result of Demengel states that any function in $\HB$ may be approximated by smooth functions.

\begin{prop}[\protect{\cite[Prop.~1.4]{demengel1984fonctions}}]
  For all $u\in \HB$, there exists a sequence $u_n\in C^\infty(\Omega)\cap W^{2,1}(\Omega)$ such that $\lim_{n\to+\infty} d(u_n,u)=0$.
  \label{prop-demengel-approx}
\end{prop}

It is possible to define the trace of a bounded Hessian function, and the normal trace of its gradient: there exists continuous (for the strong topology) linear operators $\gamma_0: \HB \rightarrow W^{1,1}(\Gamma)$ and $\gamma_1: \HB \rightarrow L^1(\Gamma)$ such that for all $u\in C^2(\overline{\Omega})$, $\gamma_0(u)=u_{|\Gamma}$ and $\gamma_{1}(u)=\Dpart{u}{\nu}_{|\Gamma}$ (where $\nu$ is the outer unit normal to $\Omega$).
It turns out that $\gamma_0$ and $\gamma_1$ are continuous w.r.t. the intermediate topology defined by~\eqref{eq-hb-distintermed}.
For the sake of simplicity, for $u\in \HB$ and when the context is clear, we shall sometimes write $u_{|\Gamma}$ (resp. $\Dpart{u}{\nu}_{|\Gamma}$) instead of $\gamma_0(u)$ (resp. $\gamma_1(u)$).
It is worth noting that one may require in Proposition~\ref{prop-demengel-approx} (see \cite{demengel1984fonctions}) that the approximating sequence $(u_n)_{n\in\NN}$ should consist of smooth functions which have the same traces as $u$: $\gamma_0 u_n = \gamma_0 u$, $\gamma_1 u_n = \gamma_1 u$ for all $n\in\NN$.

Eventually, we shall rely on the following gluing theorem proved by Demengel~\cite{demengel1984fonctions} (see also~\cite[Corollary 3.89]{Ambrosio}).
Let $U$, $V\subset \RR^2$ be two open sets with Lipschitz boundary such that $U\subset V$, and let $\Gamma=V\cap \partial U$. For $u\in \mathrm{HB}(U)$, $v\in \mathrm{HB}(V\setminus\overline{U})$ define 
\begin{align*}
  w=\left\{\begin{array}{ll} u & \mbox{ in } U,\\
    v & \mbox{ in } V\setminus \overline{U}.\end{array}\right.
\end{align*}

\begin{theorem}[\cite{demengel1984fonctions}]
 The function $w$ is in $\mathrm{HB}(V)$ if and only if $\gamma_0 u = \gamma_0 v$ on $\Gamma$, and then
  \begin{align}
    D^2 w= (D^2 {u}) \restr{U} + (D^2 {v}) \restr{(V\setminus\overline{U})} + \left(\Dpart{v}{\nu}-\Dpart{u}{\nu}\right)\nu\otimes \nu\ \Hh^1\restr \Gamma,
  \end{align}
  where $\nu$ is the normal from $U$ to $V\setminus \overline{U}$.
  \label{thm-gluing}
\end{theorem}

\subsection{Functionals Involving the Hessian}

\subsubsection{Convex Function of a Measure}
\label{sec-convex-measure}

In this paragraph, we recall some results about convex functionals of a measure. We refer the reader to~\cite[Section~2.6]{Ambrosio} for a detailed exposition on this notion.

We consider a Borel function $f: \tilde{\Omega}\times \RR^{N}\rightarrow [0,\infty)$ (where $N\in\NN^*$ and $\tilde{\Omega}$ is defined in Section~\ref{sec-domain-assumption}), positively $1$-homogeneous and convex in the second variable. Given any bounded Radon measure $\mu\in \Mm_b(\tilde{\Omega},\RR^N)$, $\mu$ is absolutely continuous w.r.t.~its total variation $|\mu|$ and the Radon-Nikodym derivative $ \frac{\mu}{|\mu|}$ is a $|\mu|$-measurable function (such that $\frac{\mu}{|\mu|}(x)\in\mathbb{S}^{N-1}$ for $|\mu|$-almost every $x\in \tilde{\Omega}$). Thus we may define
\begin{align}
  J(\mu)&\eqdef\int_{\tilde{\Omega}} f\left(x, \frac{\mu}{|\mu|}(x)\right)d|\mu|(x).
\end{align}
One may show that $J$ is convex, positively $1$-homogeneous, and that $J(\mu_1+\mu_2) = J(\mu_1)+J(\mu_2)$ when the measures $\mu_1$ and $\mu_2$ are mutually singular. Moreover, a result by Reshetnyak~\cite[Th.~2.38]{Ambrosio} ensures that it is weak-* lower semi-continuous provided $f:\tilde{\Omega}\times \RR^N\rightarrow [0,\infty]$ is lower semi-continuous. More precisely, if $\mu_n\in \Mm_b(\tilde{\Omega},\RR^N)$ and $\mu_n \stackrel{*}{\rightharpoonup} \mu\in \Mm_b(\tilde{\Omega},\RR^N)$, then 
\begin{align}
  J(\mu)\leq \liminf_{n\to +\infty} J(\mu_n). 
\end{align}
If additionally the restriction of $f$ to $\tilde{\Omega}\times \mathbb{S}^{N-1}$ is continuous and bounded, and that $|\mu_n|(\tilde{\Omega})\rightarrow |\mu|(\tilde{\Omega})$, Reshetnyak's continuity theorem (\cite[Th.~2.39]{Ambrosio}) states that
\begin{align}
  J(\mu)= \lim_{n\to +\infty} J(\mu_n). 
\end{align}

\subsubsection{Energy with Boundary Terms}
\label{sec-energy}

Now we assume that $\pi$ satisfies the hypotheses of Section~\ref{sec-pi-hypo}.
From Theorem~\ref{thm-gluing}, we know that we may extend any function $u\in \HB$ such that $\gamma_0 u=0$ in $\Gamma_D$ into a function $u\in \mathrm{HB}(\tilde{\Omega})$ such that $u=0$ in $\tilde{\Omega}\setminus \overline{\Omega}$. Combining this observation with the framework of Section~\ref{sec-convex-measure}, we define its energy as $\int_{\tilde{\Omega}} \pi(x,\frac{D^2u}{|D^2u|}(x))|D^2u|$, which turns out to be
\begin{align}
  J(u)&=\int_{\Omega} \pi\left(x,\frac{D^2u}{|D^2u|}(x)\right)|D^2u| + \int_{\Gamma_N} \pi\left(x, -\Dpart{u}{\nu}\nu\otimes \nu\right)d\Hh^1(x).
\end{align}

The next two propositions describe the (lower semi-)continuity properties of $J$ with respect to the topologies of $\HB$.

\begin{prop}
 The functional $J:\HB \rightarrow [0,\infty)$ is sequentially lower semi-continuous for the weak topology of $\HB$ (or the $L^1(\Omega)$ strong topology).
   \label{prop-J-sci}
\end{prop}

\begin{proof}
  Let $(u_n)_{n\in\NN} \in (\HB)^\NN$ which converges towards some $u\in\HB$ for the weak topology (or the $L^1(\Omega)$ strong topology). We want to prove that
  \begin{align*}
    \liminf_{n\to+\infty} J(u_n) \geq J(u).
  \end{align*}
  Recall that $J(u)=\int_{\tilde{\Omega}}\pi\left(x,\frac{D^2u}{|D^2u|}(x)\right)|D^2u|$, where $u$ is extended as the null function on $\tilde{\Omega}\setminus \overline{\Omega}$.
  However we cannot directly apply Reshetnyak's lower semi-continuity theorem since the sequence $(D^2 u_n)$ does not \textit{a priori} converge towards $D^2u$ for the weak-* topology of $\Mm_b(\tilde{\Omega})$.

Let $l=\liminf_{n\to+\infty} J(u_n)$.
If $l=+\infty$, there is nothing to prove. Assuming that $l<+\infty$, we extract a subsequence $u_{n'}$ such that $\lim_{n\to\infty} J(u_{n'})=l$.
Then for $n'$ large enough,
\begin{align}
  \alpha |D^2u_{n'}|(\tilde{\Omega}) =\alpha \left(|D^2u_{n'}|(\Omega) + \int_{\Gamma_N} \left|\Dpart{u_{n'}}{\nu} \right|d\Hh^1  \right) \leq J(u_{n'})\leq l+1.
\end{align}
Since $u_{n'}$ converges towards $u$ in $L^1(\tilde{\Omega})$, we obtain that $D^2u_n$ converges towards $D^2u$ for the weak-* topology of $\Mm_b(\tilde{\Omega})$. By Reshetnyak's lower semi-continuity theorem~\cite[Th.~2.38]{Ambrosio}, we obtain the desired inequality.
\end{proof}

\begin{prop}
 The functional $J:\HB \rightarrow [0,\infty)$ is continuous for the intermediate topology of $\HB$.
   \label{prop-J-continue}
\end{prop}

\begin{proof}
  Let $(u_n)_{n\in\NN}\in (\HB)^\NN$ be a sequence such that $\lim_{n\to \infty} d(u,u_n)=0$, where $d$ is the distance defining the intermediate topology of $\Omega$ (see~\eqref{eq-hb-distintermed}).
  
 By continuity of the trace $\gamma_1$ with respect to the intermediate topology, $\lim_{n\to+\infty}\gamma_1 u_n=\gamma_1 u$ in $L^1(\partial \Omega)$, so that 
  \begin{align*}
|D^2u_n|(\Omega)+\int_{\Gamma_N} \left|\Dpart{u_n}{\nu} \right|d\Hh^1 &\to |D^2u|(\Omega)+\int_{\Gamma_N} \left|\Dpart{u}{\nu} \right|d\Hh^1,
\end{align*}
or equivalently, considering the extension of $u$ to $\tilde{\Omega}$ described above,
\begin{align*}
  \qquad \quad |D^2u_n|(\tilde{\Omega}) &\to |D^2u|(\tilde{\Omega}).
  \end{align*}
  Moreover, since $\lim_{n\to +\infty} u_n= u$ in $L^{1}(\tilde{\Omega})$, we see that $D^2u_n\stackrel{*}{\rightharpoonup} D^2u$ in $\Mm_b(\tilde{\Omega})$.
  
 By Reshetnyak's continuity theorem~\cite[Th.~2.39]{Ambrosio}  we obtain that 
  \begin{align*}
    \lim_{n\to +\infty} \int_{\tilde{\Omega}} \pi\left(x, \frac{D^2u_n}{|D^2u_n|}\right)|D^2u_n| = \int_{\tilde{\Omega}} \pi\left(x, \frac{D^2u}{|D^2u|}\right)|D^2u|,
  \end{align*}
which is the desired result.
\end{proof}


\section{The Continuous Problem in $\HB$}
\label{sec-continuous-pb}

\subsection{Formulation of the Problem}

From now on, we consider $\Omega$, $\Gamma_D$, $\Gamma_N$ and $\pi$ which satisfy the hypotheses of Section~\ref{sec-preliminaries}, and we assume that we are given a load, that is a linear form $L:\HB \rightarrow \RR$ such that either $L$ is continuous with respect to the weak topology defined in Section~\ref{sec-hb-facts} or its restriction to bounded subsets of $\HB$ is (sequentially) continuous for the weak topology. This holds if for instance one of the three conditions hold:
\begin{itemize}
  \item $L\in (W^{1,1}(\Omega))'$,
  \item $L$ is the second derivative of a continuous function with compact support in $\Omega$, e.g. $L=\div \div M$ for some $M\in C_c(\Omega,\tens)$,
  \item $L$ is a measure which does not charge horizontal and vertical lines,
\end{itemize}
or if $L$ is a linear combination of such linear forms. 

\begin{remark}
  It is important to observe that not all measures $\mu\in\Mm_b(\Omega)$ are continuous for the $\HB$ weak topology. In~\cite{demengel1989compactness}, the author has exhibited a bounded sequence of functions $u_n$ which converges to $0$ for the $\HB$ weak topology but such that $\langle \delta_{x_0}, u_n\rangle =1$ for all $n\in\NN$ (where $\delta_{x_0}$ is the Dirac mass at some $x_0\in\Omega$).
\end{remark}

The continuity of $L$ in the first two cases is straightforward, let us give some details on the last case. 
Let $L$ be a finite measure which does not charge horizontal and vertical lines (in the sense that neither its positive nor negative part does) and $u_n$ converge weakly in $\HB$ to some $u$ and \textit{such that $\|u_n\|_{\HB}$ is bounded}. We claim then that $\langle L, u_n\rangle$ converges  to $\langle L, u\rangle$. Proceeding as in \cite{demengel1984fonctions}, by suitable extension arguments, we may assume that $u_n$ are $\mathrm{HB}$ in the whole plane with (a common) compact support. We then have (see \cite{demengel1984fonctions})  $\sup_n \Vert u_n \Vert_{L^\infty} \leq C$ and 
\[ u_n(x,y)=\theta_n(Q_{x,y}), \; u(x,y)=\theta(Q_{x,y})\]
where  $\theta_n$ and $\theta$ are the measures of mixed second partial derivatives
\[\theta_n\eqdef\frac{\partial^2u_n}{\partial x\partial y}, \;  \theta\eqdef \frac{\partial^2u}{\partial x\partial y}, \; Q_{x,y}\eqdef(-\infty, x]\times(-\infty, y].\]
By assumption $\theta_n$ converges weakly $*$ to $\theta$, but since $\vert \theta_n\vert$ is bounded, we may also assume, taking a subsequence if necessary, that $\vert \theta_n\vert$  converges weakly $*$ to some measure $\mu$ (with $\mu\ge \vert \theta\vert$). There exists two subsets of $\RR$, $I$ and $J$ which are at most countable and for which as soon as $x\in \RR\setminus I$ and $y\in \RR\setminus J$ one has:
\begin{equation}\label{bordnegl}
\mu(\partial Q_{x,y})=0.
\end{equation}
Now we claim that, as soon as $x\in \RR\setminus I$ and $y\in \RR\setminus J$,  $u_n(x,y)$ converges to $u(x,y)$. Indeed, let $\varepsilon>0$ and $\phi_\varepsilon\in C(\RR^2, [0,1])$ vanishing outside $Q_{x,y}$ and equal to $1$ on $Q^{\varepsilon}_{x,y}:=(-\infty, x-\varepsilon)\times (-\infty, y-\varepsilon)$, we then have for every $\varepsilon >0$:
\[\begin{split}
\limsup_n \vert u_n(x,y)-u(x,y)\vert \le \limsup_n \vert \langle \theta_n -\theta, \phi_ \varepsilon  \rangle \vert \\
+ \limsup_n \vert \theta_n \vert (Q_{x,y}\setminus Q_{x,y}^\varepsilon)+ \vert \theta \vert   (Q_{x,y}\setminus Q_{x,y}^\varepsilon)\\
\le 2 \mu( Q_{x,y}\setminus Q_{x,y}^\varepsilon)
\end{split}\]
where the last line follows from the weak $*$ convergence of $\theta_n$ to $\theta$, the fact that  $(Q_{x,y}\setminus Q_{x,y}^\varepsilon)$ is closed and the inequality $\mu\ge \vert \theta\vert$. We then let $\varepsilon$ tend to $0$ and use \eqref{bordnegl} to deduce that $u_n(x,y)$ converges to $u(x,y)$. This proves that $u_n$ converges to $u$ $\vert L\vert$-a.e., and together with the $L^{\infty}$ bound and Lebesgue's dominated convergence theorem, this yields the claimed convergence.

\begin{remark}
Because of the coercivity of $J$, all the sequences considered in this paper shall be bounded in~$\HB$. Therefore, in the following, we shall often refer to the ``continuity of $L$ for the weak topology'', tacitly meaning that the given argument also holds for loads $L$ such that their restrictions to bounded subsets of $\HB$ are continuous for the weak topology.
\end{remark}

We also assume that there exists a function $\hat{u}\in\HB$ such that $\hat{u}_{\vert \Gamma_D}=0$ and $\langle L,\hat{u}\rangle \neq 0$.


The problem we want to solve is: 
\begin{align}
  \inf_{u\in \HB} J(u) \quad &\mbox{such that } \left\{\begin{array}{l} u=0 \mbox{ on } \Gamma_D\\
    \langle L, u\rangle =1   \end{array}\right.\tag{$\LA$}\label{eq-prelim-problem}
\\
\mbox{ where }\quad J(u)&= \int_{\Omega} \pi\left(x,\frac{D^2u}{|D^2u|}(x)\right)|D^2u| +\int_{\Gamma_N} \pi\left(x,-\frac{\partial u}{\partial \nu}\nu\otimes \nu\right)d\Hh^1.\nonumber
\end{align}

We may also write~\eqref{eq-prelim-problem} as the problem of minimizing the functional $F$ over $\HB$, where
\begin{align}
   F(u)&\eqdef\left\{\begin{array}{ll} J(u) &\mbox{ if } u_{\Gamma_D}=0 \mbox{ and } \langle L,u\rangle =1,\\
    +\infty \mbox{ otherwise.}\end{array}\right.
\end{align}

\subsection{Existence of a Minimizer}

We obtain the existence of a minimizer to~\eqref{eq-prelim-problem} by the direct method of the calculus of variations.
\begin{prop}[Existence of a minimizer]
There exists a solution $u^\star\in \HB$ to Problem~\eqref{eq-prelim-problem}.
\end{prop}

\begin{proof}
  Let $(u_n)_{n\in \NN}\in \HB^\NN$ be any minimizing sequence. Problem~\eqref{eq-prelim-problem} is feasible since there exists $\hat{u}$ such that $\hat{u}_{|\Gamma_D}=0$ and $\langle L,\hat{u}\rangle \neq 0$. Hence there exists $M\in (0,+\infty)$ such that for all $n\in \NN$, $0\leq J(u_n)\leq M$. This implies by the coercivity of the functional (see Lemma~\ref{lem-poincare} below) that $\normhb{u_n}$ is bounded. 

  Therefore, we may extract a subsequence $(u_{n'})_{n'\in\NN}$ which converges to some $u^\star\in \HB$ for the weak topology.
 From the continuity of the trace operator $\gamma_0: W^{1,1}(\Omega)\rightarrow L^1(\partial \Omega)$, we obtain that $u^\star= 0$ on $\Gamma_D$.
 Moreover, from the continuity of $L$ for the weak topology, $\langle L, u^\star\rangle =1$.  By Proposition~\ref{prop-J-sci} we obtain $J(u^\star)\leq \liminf_{n'\to+\infty} J(u_{n'})=\inf_{u\in\HB}F(u)$.

 Hence $u^\star$ is a solution to Problem~\eqref{eq-prelim-problem}.
\end{proof}

In the above proof, we have used the coercivity of $J$ with respect to the $\HB$ norm. This is a consequence of the following Lemma.
\begin{lem}[Coercivity]
There exists $C>0$ (which depends only on $\Omega$, $\pi$ and $\Gamma_D$) such that for all $u\in \HB$ with $u=0$ on $\Gamma_D$:
\begin{align}
  J(u) \geq C\normhb{u}=C\left(\|u\|_{W^{1,1}}+ |D^2u|(\Omega) \right).
\label{eq-basic-poincare}
\end{align}
\label{lem-poincare}
\end{lem}

\begin{proof}
  We follow the standard proof of the Poincar\'e inequality (see \cite{Ambrosio}).
Assume by contradiction that \eqref{eq-basic-poincare} does not hold. Then there exists a sequence $(u_n)_{n\in \NN} \in (\HB)^\NN$, such that for all $n\in \NN$, $\gamma_0 u_n =0$ on  $\Gamma_D$ and
\begin{align}	
  \alpha \int_{\Omega} |D^2u_n| \leq \int_\Omega \pi\left(x,\frac{D^2u_n}{|D^2u_n|}\right)|D^2u_n|  < \frac{1}{n} \|u_n\|_{\HB}.
\end{align}
Since $\normhb{u_n}\neq 0$, we may assume, up to a rescaling, that $\normhb{u_n}=1$ so that the sequence is bounded in 
$\HB$.

Hence (see Section~\ref{sec-hb-facts}), we may extract a subsequence $(u_{n'})_{n'\in\NN}$ which weakly converges to some $\bar{u}\in \HB$. Since weak convergence in $\HB$ implies strong convergence in $W^{1,1}(\Omega)$ we have $\gamma_0 \bar{u}= 0$ on $\Gamma_D$.
 By the lower semi-continuity of the second variation, 
\begin{align*}
  0\leq \alpha\int_\Omega |D^2\bar{u}| \leq \liminf_{n'\to +\infty} \alpha\int_\Omega |D^2u_{n'}| =0,
\end{align*}
so that $\bar{u}\in \Pp_1$. Since $\Gamma_D$ contains three points which are not aligned, this implies that $\bar{u}=0$, which contradicts $\lim_{n'\to +\infty} \normhb{u_{n'}}=1$, 
hence the claimed result. 
\end{proof}

\subsection{An Approximation Result}

We shall rely on the following result in order to prove the $\Gamma$-convergence of the finite element method.

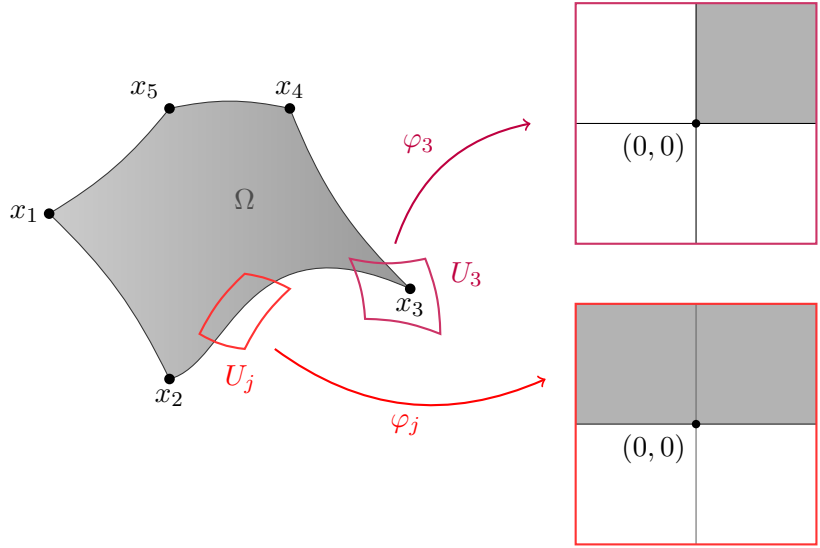
\begin{figure}
  \centering
  \begin{tikzpicture}[scale=2.0]
\definecolor{vert}{rgb}{0.15 , 0.5 , 0.25} 
\definecolor{vert}{rgb}{0.1725 , 0.6627 , 0.1725} 
 \begin{scope}[shift={(0,0)}]
   \coordinate (A) at (-0.3,0.4) ;
  \coordinate (B) at (0.5,-0.7) ;
  \coordinate (C) at (2.1,-0.1) ;
  \coordinate (D) at (1.3,1.1);
  \coordinate (E) at (0.5,1.1);
  \coordinate (O) at (1.0,0.5);
  \coordinate (Abis) at (0.2,0.7);
  \coordinate (Bbis) at (0.9,-0.6);
  \coordinate (Bter) at (1.0,0.4);
  \path[left color=gray!40, right color=gray!80,draw=black!80] (A) node[left] {$x_1$} to[bend left=10] (B) node[below] {$x_2$} ..controls (Bbis) and  (Bter)..  (C) node[below] {$x_3$} to[bend left=12] (D) node[above] {$x_4$} to[bend right=12] (E) node[above left] {$x_5$} to[bend left=10] (A) -- cycle;
  \draw[color=black!70] (O) node {$\displaystyle{\Omega}$};
  \fill (A)  circle[radius=1pt];
  \fill (B)  circle[radius=1pt];
  \fill (C)  circle[radius=1pt];
  \fill (D)  circle[radius=1pt];
  \fill (E)  circle[radius=1pt]; 
  \coordinate (a) at (2.3,-0.4);
  \coordinate (b) at (2.2,0.1);
  \coordinate (c) at (1.7,0.1);
  \coordinate (d) at (1.8,-0.3);
  \path[draw=purple!80,thick] (a) to[bend right=10] (b) to[bend left=12] (c) to [bend left=10] (d) to[bend left=10] (a) --cycle;
  \coordinate (X) at (2.0,0.2);
  \coordinate (Z) at (2.9,1.0);
  \draw [->,color=purple,thick] (X) to [bend left=30] node[midway,above left] {$\varphi_{3}$} (Z);
  \draw[color=purple] (2.3,0.0) node[right] {$U_3$} ;
  \coordinate (a) at (1.3,-0.1);
  \coordinate (b) at (1.0,-0.5);
  \coordinate (c) at (0.7,-0.4);
  \coordinate (d) at (1.0,-0.0);
  \path[draw=red!80,thick] (a) to[bend right=10] (b) to[bend left=12] (c) to [bend left=10] (d) to[bend left=10] (a) --cycle;
  \coordinate (X) at (1.2,-0.5);
  \coordinate (Z) at (3.0,-0.7);
  \draw [->,color=red, thick] (X) to [bend right=30] node[midway,below] {$\varphi_j$} (Z);
  \draw[color=red] (0.8,-0.7) node[right] {$U_j$} ;
  \end{scope}

 \begin{scope}[shift={(4,1)},scale=0.8]
   \coordinate (A) at (-1,-1) ;
   \coordinate (B) at (1,-1) ;
   \coordinate (C) at (1,1) ;
   \coordinate (D) at (-1,1) ;
   \coordinate (E) at (-1,0) ;
   \coordinate (F) at (1,0) ;
   \coordinate (G) at (0,-1) ;
   \coordinate (H) at (0,1) ;
   \coordinate (O) at (0,0) ;
   \path[fill,color=gray!60] (H) to (C) to (F) to (O) to (H) -- cycle;
   \draw (E) -- (F);
   \draw (G) -- (H);
   \draw (O) node[below left] {$(0,0)$};
   \fill (O)  circle[radius=1pt]; 
   \path[draw=purple!80,thick] (A) to (B) to (C) to (D) to (A) -- cycle;
  \end{scope}
  
  \begin{scope}[shift={(4,-1)},scale=0.8]
   \coordinate (A) at (-1,-1) ;
   \coordinate (B) at (1,-1) ;
   \coordinate (C) at (1,1) ;
   \coordinate (D) at (-1,1) ;
   \coordinate (E) at (-1,0) ;
   \coordinate (F) at (1,0) ;
   \coordinate (G) at (0,-1) ;
   \coordinate (H) at (0,1) ;
   \coordinate (O) at (0,0) ;
   \path[fill=gray!60,draw=black!80] (C) to (D) to (E) to[thick] (F) to (C) -- cycle;
   \draw[color=black!60] (E) -- (F);
   \draw[color=black!60] (G) -- (H);
   \draw (O) node[below left] {$(0,0)$};
   \fill (O)  circle[radius=1pt]; 
   \path[draw=red!80,thick] (A) to (B) to (C) to (D) to (A) -- cycle;
  \end{scope}
\end{tikzpicture}
  \caption{The local charts used to smooth the function $u$ in Proposition~\ref{prop-approx-hb}.
}
  \label{fig-chart}
\end{figure}

\begin{prop}
  Let $u\in \HB$ and $\epsilon>0$, there exists $u_\epsilon\in C^3(\overline{\Omega})$ such that
  \begin{align}
    u_{\Gamma_D}=0,\\
    \|u-u_{\epsilon}\|_{W^{1,1}(\Omega)}\leq \epsilon,\\
    \left| |D^2u|(\Omega)-|D^2u_{\epsilon}|(\Omega)\right|\leq \epsilon,\\
    \left| J(u)-J(u_\epsilon)\right| \leq \epsilon.
  \end{align}
  \label{prop-approx-hb}
\end{prop}

\begin{proof}
  From Propositions~\ref{prop-demengel-approx} and~\ref{prop-J-continue}, we see that we can already find a function $v\in W^{2,1}(\Omega)$ such that $v_{\vert \Gamma_D}=0$, $\|u-v\|_{W^{1,1}}\leq \frac{\varepsilon}{2}$, $\left| |D^2u|(\Omega)-|D^2v|(\Omega)\right|\leq \frac{\epsilon}{2}$, and $|J(u)-J(v)|\leq \frac{\varepsilon}{2}$. Therefore, replacing $u$ with $v$, there is no loss of generality in assuming that $u\in W^{2,1}(\Omega)$.

  \paragraph{Partition of the domain.}
  By assumption on $\Omega$ (Section~\ref{sec-domain-assumption}), there exists a finite open cover $(U_i)_{0\leq i\leq N}$ of $\overline{\Omega}\subset \RR^2$, such that 
  $U_0\subcomp \Omega$, and for each $i\in \{1, \ldots, N\}$ there exists a $C^3$-diffeomorphism $\varphi_i: U_i \rightarrow Q$ (where $Q=(-1,1)\times (-1,1)$) which satisfies (see Figure~\ref{fig-chart})
  \begin{itemize}
    \item for $1\leq i \leq J$,  $\varphi_i(\overline{\Omega}\cap U_i)=[0,1)\times [0,1)$, where $\varphi_i^{-1}((0,0))=x_i$ and $x_i$ is one of the exceptional points where $\partial \Omega$ is not smooth (see Section~\ref{sec-domain-assumption}),
  \item for $J+1\leq i\leq N$,  $\varphi_i(\overline{\Omega}\cap U_i)= (-1,1)\times [0,1)$.
\end{itemize}

Let us choose a $C^\infty$-partition of unity $(\beta_i)_{0\leq i\leq N}$ of $\overline{\Omega}$ with respect to the cover $(U_i)_{0\leq i\leq N}$.
Each function $u_i\eqdef  \beta_i u$ is in $W^{2,1}(\Omega)$ with support in $\supp \beta_i\cap \overline{\Omega }$, and $u_i$ satisfies the boundary condition $u_i=0$ on $\Gamma_D\cap U_i$. 

Now we define an extension and smoothing for each $u_i$ depending on the required boundary condition: we must not introduce discontinuities at $\partial \Omega\cap U_i$, but we must also preserve the property that that $u_i=0$ on $\Gamma_D\cap U_i$.

For $J+1\leq i\leq N$, $U_i$ intersects exactly one connected component  of $\partial \Omega\setminus \{x_1,\ldots x_J\}$, so that $U_i\cap \Gamma_D=\emptyset$, or $U_i\cap \Gamma_D\neq\emptyset$ and $U_i\cap \Gamma_D=U_i\cap \partial \Omega$. 
Similarly, 
 each $U_i$ for $1\leq i\leq J$ intersects (at most) two connected components of $\partial \Omega\setminus \{x_1,\ldots x_J\}$, each of which may be involved in the Dirichlet condition or not.
For the sake of simplicity, we detail the argument in the case $J+1\leq i\leq N$. The case $1\leq i\leq J$ is handled by applying the argument below in both the horizontal and vertical axes.

\paragraph{Extension.}

Let us write $Q^+=(-1,1)\times (0,1)$ and set $v_i=u_i \circ \phi_i^{-1} : Q^+\rightarrow \RR$. If $\Gamma_D\cap U_i=\emptyset$, we define, for $(x,t)\in Q=(-1,1)\times (-1,1)$
\begin{align}
   \hat{v}_i(x,t) &=\left\{\begin{array}{l} 
    v_i(x,t) \mbox{ for } t>0,\\
    -2v_i(x,-2t)+3v_i(x,-t) \mbox{ for } t<0.
\end{array}\right.
\label{eq-continuous-extension1}
\end{align}
In the case where $\Gamma_D\cap U_i=\partial \Omega\cap U_i$, we use a different extension:
\begin{align}
   \hat{v}_i(x,t) &=\left\{\begin{array}{l} 
    v_i(x,t) \mbox{ for } t>0,\\
    -v_i(x,-t) \mbox{ for } t<0.
\end{array}\right.
\label{eq-continuous-extension2}
\end{align}
Since $v_i\in W^{2,1}\left(Q^+\right)$, Theorem~\ref{thm-gluing} ensures in both cases that $\hat{v}_i \in W^{2,1}(Q)$ (since the normal trace of the gradient has no jump along $(-1,1)\times \{0\}$) with compact support in $Q$.

\paragraph{Smoothing.} 

Now we define an approximation to the identity $(\rho_h)_{h>0}$ on $Q$. Given some even function $\eta \in C_c^\infty(\RR,[0,+\infty))$ such that $\supp \eta \subset (-1,1)$ and $\int_\RR \eta =1$, we set $\rho:\begin{array}{rccc} &\RR^2&\longrightarrow &\RR_+\\ &(y_1,y_2)&\longmapsto &\eta(y_1)\eta(y_2)\end{array}$ and $\rho_h= \frac{1}{h^2}\rho(\frac{\cdot}{h})$ for all $h>0$. 
  We set $\tilde{v}_i=\rho_{h_i}\ast \hat{v}_i$ for $h_i<\frac{1}{2} \dist(\supp \hat{v_i}, \partial Q)$, so that $\tilde{v}_i\in C^{\infty}_c(Q)$ and 
  \begin{align*}
    \lim_{h_i\to 0^+} \|v_i- \tilde{v}_i\|_{W^{2,1}(Q^+)}=0.
  \end{align*}
  Moreover, we observe that in the case where $\Gamma_D\cap U_i=\partial \Omega\cap U_i$, Equation~\eqref{eq-continuous-extension2} implies that $\tilde{v}_i=0$ in $(-1,1)\times \{0\}$.

\paragraph{Back to $U_i$ and $\Omega$.} 

It is easily seen that $\tilde{u}_i\eqdef \tilde{v}_i\circ \varphi_i\in C_c^3(U_i)$, and $\tilde{u}_i=0$ on $\Gamma_D\cap U_i$. Moreover there exists a constant $C>0$ which depends only on $\varphi_i$ such that 
\begin{align*}
  \| \tilde{u}_i - u_i\|_{W^{2,1}(U_i\cap \Omega)}\leq C \|\tilde{v}_i-v_i\|_{W^{2,1}(Q^+)},
\end{align*}
hence $\| \tilde{u}_i - u_i\|_{W^{2,1}(U_i\cap \Omega)} \to 0$ as $h_i\to 0$. As a consequence, we may choose $(h_i)_{1\leq i \leq N}$ such that
\begin{align*}
  \left\| u-\sum_{i=1}^N \tilde{u}_i\right\|_{W^{2,1}(\Omega)}\leq \sum_{i=1}^N \| u_i-\tilde{u}_i\|_{W^{2,1}(\Omega)}
\end{align*}
is arbitrarily small. The function $\sum_{i=1}^N \tilde{u}_i$ is in $C^3(\overline{\Omega})$, and by Proposition~\ref{prop-J-continue}, we obtain the claimed inequalities.
\end{proof}


\section{Finite Elements Discretization}
\label{sec-elfini}

This section is devoted to the analysis of the finite element method for Problem~\eqref{eq-prelim-problem}. We refer to~\cite{brenner2008} for a comprehensive exposition of the theory of finite elements.

\subsection{Notations and Definitions}
\label{sec-fem-notations}

From now on, we assume that $ \Omega\subset \RR^2$ is a bounded convex polyhedral open set (which implies the assumptions of Section~\ref{sec-domain-assumption}.
We say that $\Tt$ is a triangulation of $\Omega$, if it is a finite collection of open triangles $\{T_i\}$ such that
\begin{itemize}
  \item $T_i\cap T_j=\emptyset$ for $i\neq j$,
  \item $\bigcup \overline{T_i} = \overline{\Omega}$,
\end{itemize}
and no vertex of any triangle lies in the interior of an edge of another triangle.

We consider a family of triangulations of $\Omega$, $\{\Ttn, n\in\NN\}$ such that there exists a decreasing sequence $(h_n)_{n\in\NN}\in\RR^\NN$, such that $\lim_{n\to +\infty}h_n =0$ and
\begin{align}
  \max\enscond{\diam T}{T\in\Ttn}\leq  h_n\diam \Omega. \label{eq-h-def}
\end{align}
We say that the family is \textit{nondegenerate} if there exists $\rho>0$ such that for all $n\in\NN$ and $T\in \Ttn$,
\begin{align}
  \diam B_T&\geq \rho \diam T\label{eq-rho-def}
\end{align}
where $B_T$ is the largest ball contained in $T$.

Given a triangulation $\Tt=\Ttn$ for some $n\in\NN$ and $k\in\NN$, the finite element method consists in defining a global interpolant $\IT: C^k(\overline{\Omega})\rightarrow \RR^\Omega$ which coincides with some local interpolant $\It$ on each triangle $T\in \Tt$: $({\IT}u)_{\vert T} \eqdef \It u$. We describe below the finite elements and the corresponding local interpolants used in this paper.

\paragraph{Lagrange finite elements.}

\begin{figure}
  \centering
  \begin{tikzpicture}
 \begin{scope}[shift={(0,0)}]
   \coordinate (A) at (0,0) ;
  \fill (A)  circle[radius=2pt];
  \coordinate (B) at (2,0) ;
  \fill (B)  circle[radius=2pt];
  \coordinate (C) at (0,2) ;
  \fill (C)  circle[radius=2pt];
  \coordinate (D) at ($0.5*(A)+0.5*(B)$);
  \fill (D)  circle[radius=2pt];
  \coordinate (E) at ($0.5*(B)+0.5*(C)$);
  \fill (E)  circle[radius=2pt];
  \coordinate (F) at ($0.5*(C)+0.5*(A)$);
  \fill (F)  circle[radius=2pt]; 
  \draw (A) node[left] {$z_{0,1}$} -- (B) node[right] {$z_{0,2}$} -- (C) node[above] {$z_{0,3}$} --cycle ;
  \draw (D) node[below] {$z_{0,4}$}  (E) node[right] {$z_{0,5}$}  (F) node[left] {$z_{0,6}$}  ;
  \coordinate (T) at (1.5,2);
  \draw (T) node {$T_0$};
  \end{scope}
 \begin{scope}[shift={(6,0)}]
   \coordinate (A) at (0,0) ;
  \fill (A)  circle[radius=2pt];
  \coordinate (B) at (2.4,0.2) ;
  \fill (B)  circle[radius=2pt];
  \coordinate (C) at (-0.6,1.2) ;
  \fill (C)  circle[radius=2pt];
  \coordinate (D) at ($0.5*(A)+0.5*(B)$);
  \fill (D)  circle[radius=2pt];
  \coordinate (E) at ($0.5*(B)+0.5*(C)$);
  \fill (E)  circle[radius=2pt];
  \coordinate (F) at ($0.5*(C)+0.5*(A)$);
  \fill (F)  circle[radius=2pt]; 
  \draw (A) node[left] {$z_1$} -- (B) node[right] {$z_2$} -- (C) node[above] {$z_3$} --cycle ;
  \draw (D) node[below] {$z_4$}  (E) node[above right] {$z_5$}  (F) node[left] {$z_6$}  ;
  \coordinate (T) at (0.6,1.7);
  \draw (T) node {$T$};
  \end{scope}

  \coordinate (X) at (2.2,1.5);
  \coordinate (F) at (3.9,2.2);
  \draw (F) node {$A$};
  \coordinate (Z) at (4.7,1.5);
\draw [->] (X) to [bend left=30] (Z);

\end{tikzpicture}
  \caption{The Lagrange $\Pp_2$ interpolation is determined by 6 control nodes which impose the values of the polynomial at $z_1,\ldots z_6$. Those nodes are the three vertices and the middle of each edge. Each element $T$ is the image of the reference element $T_0$ through some affine map $A$.}
  \label{fig-lagrange}
\end{figure}
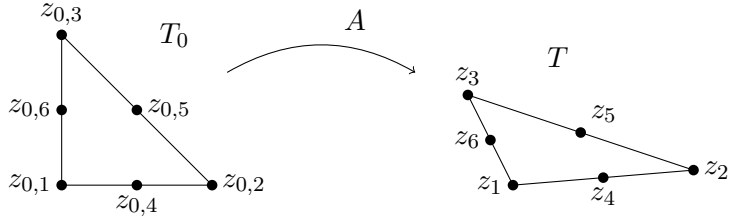

The finite elements are $(T,\Pp_2,\Nn_T)$, where $T$ is any triangle of $\Tt$, $\Pp_2$ is the space of quadratic polynomials on $\Tt$ ($\dim \Pp_2=6$) and $\Nn= \{N_1,\ldots N_{6} \}\subset C^0(\overline{T})'$ where for $1\leq i\leq 6$, $N_i$ is the evaluation at point $z_i$ (see Figure~\ref{fig-lagrange}). Observe that $\Nn$ is a basis of $\Pp_2'$.

The local interpolant $\It:C^0(\overline{T})\rightarrow \Pp_2$ is defined by 
\begin{align}
  \It(u)\eqdef \sum_{i=1}^6 N_i(u)\psi_i = \sum_{i=1}^6 u(z_i)\psi_i
\end{align}
where the collection of shape functions $(\psi_i)_{1\leq i\leq 6}$ is the basis (of $\Pp_2$) dual to $\Nn$. 

The finite elements are affine-equivalent to some reference element denoted $(T_0, \Pp_2, \Nn_{T_0})$, where for instance the vertices are $(0,0)$, $(1,0)$, and $(0,1)$ (see the definition in~\cite[Section~3.4]{brenner2008} or the Hermite case below).

\paragraph{Hermite finite elements.}

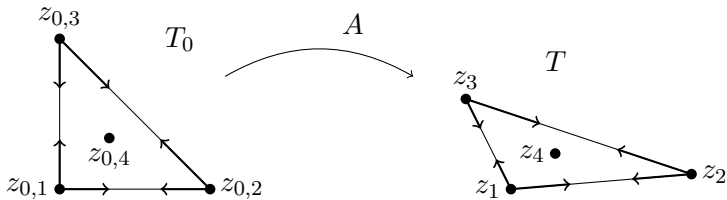
\begin{figure}
  \centering
  \begin{tikzpicture}

 \begin{scope}[shift={(0,0)}]
   \coordinate (A) at (0,0) ;
  \fill (A)  circle[radius=2pt];
  \coordinate (B) at (2,0) ;
  \fill (B)  circle[radius=2pt];
  \coordinate (C) at (0,2) ;
  \fill (C)  circle[radius=2pt];
  \coordinate (D) at ($0.33*(A)+0.33*(B)+0.34*(C)$);
  \fill (D)  circle[radius=2pt];
  \draw (A) node[left] {$z_{0,1}$} -- (B) node[right] {$z_{0,2}$} -- (C) node[above] {$z_{0,3}$} --cycle ;
  \draw (D) node[below] {$z_{0,4}$};
  \coordinate (T) at (1.6,2);
  \draw (T) node {$T_0$};
  \draw [->,thick] (A) to ($(A)!.33!(B)$);
  \draw [->,thick] (A) to ($(A)!.33!(C)$);
  \draw [->,thick] (B) to ($(B)!.33!(A)$);
  \draw [->,thick] (B) to ($(B)!.33!(C)$);
  \draw [->,thick] (C) to ($(C)!.33!(A)$);
  \draw [->,thick] (C) to ($(C)!.33!(B)$);
  \end{scope}

 \begin{scope}[shift={(6,0)}]
   \coordinate (A) at (0,0) ;
  \fill (A)  circle[radius=2pt];
  \coordinate (B) at (2.4,0.2) ;
  \fill (B)  circle[radius=2pt];
  \coordinate (C) at (-0.6,1.2) ;
  \fill (C)  circle[radius=2pt];
  \coordinate (D) at ($0.33*(A)+0.33*(B)+0.34*(C)$);
  \fill (D)  circle[radius=2pt];
  \draw (A) node[left] {$z_1$} -- (B) node[right] {$z_2$} -- (C) node[above] {$z_3$} --cycle ;
  \draw (D) node[left] {$z_4$};
  \draw [->,thick] (A) to ($(A)!.33!(B)$);
  \draw [->,thick] (A) to ($(A)!.33!(C)$);
  \draw [->,thick] (B) to ($(B)!.33!(A)$);
  \draw [->,thick] (B) to ($(B)!.33!(C)$);
  \draw [->,thick] (C) to ($(C)!.33!(A)$);
  \draw [->,thick] (C) to ($(C)!.33!(B)$);
  \coordinate (T) at (0.6,1.7);
  \draw (T) node {$T$};
  \end{scope}

  \coordinate (X) at (2.2,1.5);
  \coordinate (F) at (3.9,2.2);
  \draw (F) node {$A$};
  \coordinate (Z) at (4.7,1.5);
\draw [->] (X) to [bend left=30] (Z);

\end{tikzpicture}
  \caption{The Hermite $\Pp_3$ interpolation is determined by 4 control nodes which impose the values of the polynomial at $z_1,\ldots z_4$ and the gradient at $z_1,z_2,z_3$. Those nodes are the three vertices and the barycenter. Each element $T$ is the image of the reference element $T_0$ through some affine map $A$.}
  \label{fig-hermite}
\end{figure}

Here, the finite elements are $(T,\Pp_3,\Nn_T)$, where $T$ is any triangle of $\Tt$, $\Pp_3$ is the space of cubic polynomials on $\Tt$ ($\dim \Pp_3=10$) and $\Nn= \bigcup_{i=1}^4\{N_i\}\cup \bigcup_{1\leq i,j\leq 3, i\neq j} \{N_{i\to j}\}$, where $N_i\in C^0(\overline{T})'$ is the evaluation at point $z_i$, and $N_{i\to j}\in C^1 (\overline{T})'$ is defined as the directional derivative of $u$ along the direction of an edge $N_{i\to j}(u)=(\nabla u(z_i))(z_j-z_i)$  (see Figure~\ref{fig-hermite}). Observe that $\Nn$ is a basis of $\Pp_3'$.

The local interpolant $\It:C^0(\overline{T})\rightarrow \Pp_3$ is defined by 
\begin{align}
  \It(u)&\eqdef \sum_{i=1}^{4} N_i(u)\psi_i+ \sum_{1\leq i,j\leq 3, i\neq j} N_{i\to j}(u)\psi_{i\to j}\\
  &= \sum_{i=1}^{4} u(z_i)\psi_i+ \sum_{1\leq i,j\leq 3, i\neq j} (\nabla u(z_i)\cdot(z_j-z_i))\psi_{i\to j}
\end{align}
where $\{\psi_i\}\cup \{\psi_{i\to j}\}$ forms the basis (of $\Pp_3$) dual to $\Nn$.

Each finite element is affine-equivalent to some reference element denoted $(T_0,\Pp_3,\Nn_{T_0})$, where for instance the vertices are $(0,0)$, $(1,0)$, and $(0,1)$: more precisely there is some affine map $A: x\mapsto ax+b$ where $a\in \GL$, $b\in \RR^2$ such that:
\begin{itemize}
  \item $A(T_0)=T$
  \item $A^*\Pp_3 =\Pp_3$, where for all $f\in\Pp_3$, $A^*f=f\circ A$,
  \item $A_*\Nn_{T_0}=\Nn_{T}$ where  $A_* N(f_0)=N(A^*f_0)$ for $f_0\in\Pp_3$.
\end{itemize}

The nondegeneracy of the triangulation $\Tt$ allows to bound the distance between $\IT u$ and $u$ for smooth functions $u$ (using the fact that the affine maps $A$ for each triangle are not too ill-conditioned): we refer to~\cite{brenner2008} for more details (see also the next section).

\subsection{Approximation of a Function with Lagrange or Hermite Finite Elements}

In this section, we prove that one may approximate (for the strong $\HB$ topology) a smooth function using the Lagrange or Hermite finite elements described above. As mentioned above, we assume that the family $\{\Ttn, n\in \NN\}$ is nondegenerate.

\begin{prop}
  Let $u\in C^3(\overline{\Omega}, \RR)$. Then, there exists a constant $C\geq 0$ (which only depends on $\rho$ and $\Omega$), such that for all $n\in\NN$ and $\Tt\eqdef\Ttn$ .
\begin{align}
  \|u-\IT u\|_{W^{1,1}(\Omega)}&\leq C |\Omega| h^2 \|\nabla^3u\|_{L^\infty(\Omega)},\label{eq-approx-deriv1}\\
  |D^2u-D^2\IT u|(\Omega)&\leq C |\Omega| h \|\nabla^3u\|_{L^\infty(\Omega)},\label{eq-approx-deriv2}
\end{align}
where $h\eqdef h_n$ is the real number given in~\eqref{eq-h-def}.
\label{prop-cv-interpolation}
\end{prop}

\begin{proof}
  In the following, $C$ is a positive constant which may change from one line to another.
  We apply~\cite[Corollary~4.4.7]{brenner2008} with $m=3$, $l=2$ (though in fact our nodal variables only depend on zero and first-order derivatives) and $p=+\infty$:  
for $0\leq i\leq 2$ and for each triangle $T\subset \RR^2$, there exists constants $C^i_{\gamma, \delta}>0$  such that ,
\begin{equation}
\|\nabla^i u - \nabla^i \It u\|_{L^\infty(T)} \leq C^i_{\gamma, \sigma} (\diam T)^{3-i} \|\nabla^3 u \|_{L^\infty(T)}.
\end{equation}

The regularity constant $C^i_{\gamma, \delta}$ actually only depends on two characteristics of the triangle defined in \cite{brenner2008}: its chunkiness $\gamma = \gamma(T)$, and a Lebesgue constant $\sigma = \sigma(T)$. These two quantities depend continuously on the triangle shape, and are scaling and translation invariant. 

Given a triangle $T$, consider $\hat T= \frac{1}{\diam T}T +b_T$ with $b_T\in \RR^2$ chosen so that $\hat T$ is centered at the origin. The set of all triangles $\hat T$ with diameter $1$ such that \eqref{eq-rho-def} holds and which are centered at the origin being compact (see the proof of Theorem~4.4.20 in~\cite{brenner2008}),  we obtain that $\gamma(\hat T)$ ($=\gamma(T)$) and $\sigma(\hat T)$ ($=\sigma(T)$) are bounded. Hence $C^i_{\gamma, \sigma}$ is bounded as well, which implies that there exists a uniform constant $C>0$ such that for all $n\in \NN$ and  $T\in \Ttn$,
\begin{align}
\|\nabla^i u - \nabla^i \It u\|_{L^\infty(T)} \leq C (\diam T)^{3-i} \|\nabla^3 u \|_{L^\infty(T)}.\label{eq-proof-derivi}
\end{align}

Now, let us prove~\eqref{eq-approx-deriv2}. We have
\begin{align*}
  |D^2u-D^2\IT u|(\Omega)= \sum_{T\in \mathcal{T}}\int_T |\nabla^2 u(x)-\nabla^2\IT u(x)| dx + \sum_{e\in \Ee} \int_e \left|\llbracket \nabla\IT u \cdot \nu\rrbracket\right| d\Hh^1.
\end{align*}
On the one hand, using~\eqref{eq-proof-derivi},
\begin{align*}
\int_T |\nabla^2u(x)-\nabla^2\IT u(x)| dx &\leq |T| \|\nabla^2u-\nabla^2\It u\|_{L^\infty(T)}\\
&\leq |T|C(\diam T)\|\nabla^3u\|_{L^\infty(T)}\\
&\leq |T|C h\|\nabla^3u\|_{L^\infty(\Omega)}.
\end{align*}
On the other hand, if one edge is shared between triangles $S$ and $T$:
\begin{align*}
\int_e \left|\llbracket \nabla\IT u \cdot \nu\rrbracket\right| d\Hh^1 &\leq |e| \|\nabla (\Is u-\It u)\|_{L^\infty(e)}\\
&\leq |e|\left( \|\nabla (u-\It u)\|_{L^\infty(e)} + \|\nabla (u-\Is u)\|_{L^\infty(e)}\right).
\end{align*}
By~\eqref{eq-proof-derivi}, since $|e|\leq h\diam\Omega$ and $(\diam T)^2\leq \frac{1}{\rho^2}\left(\frac{\diam B_T}{\diam \Omega}\right)^2\leq C'|T|$:
\begin{align*}
|e| \|\nabla (u-\It u)\|_{L^\infty(e)} &\leq |e|C (\diam T)^2 \|\nabla^3u\|_{L^\infty(T)}\\
&\leq  hC|T| \|\nabla^3u\|_{L^\infty(\Omega)}.
\end{align*}
Summing the contributions of triangles and their edges, we obtain~\eqref{eq-approx-deriv2}. Equation~\eqref{eq-approx-deriv1} is obtained in a similar but simpler way, since there is no contribution of edges.
\end{proof}

\subsection{$\Gamma$-Convergence of the Finite Element Approximation}

Let $(\Ttn)_{n\in\NN}$ be a nondegenerate family of subdivisions of $\Omega$ and  $\ITn$ the corresponding interpolation operator (see Section~\ref{sec-fem-notations}).
We consider a family of discrete problems defined by restricting \eqref{eq-prelim-problem} to the elements of $\Im \ITn$:
\begin{align}
\inf_{u\in \HB} J(u)  \quad \mbox{such that } \left\{\begin{array}{l} u=0 \mbox{ on } \Gamma_D\\
    \langle L, u\rangle =1,\\
    u\in \Im \ITn.
  \end{array}\right.\tag{$\LA_n$}
\label{eq-pbm-discret}
\end{align}
We may also write~\eqref{eq-pbm-discret} as the minimization of $F_n$, where for all $u\in \HB$,
\begin{align}
    F_n(u)&\eqdef\left\{\begin{array}{ll} J(u) &\mbox{ if } u_{\Gamma_D}=0, \langle L,u\rangle =1 \mbox{ and } u\in \Im \ITn,\\
      +\infty \mbox{ otherwise.}\end{array}\right.
\end{align}
We assume that the triangulation $\Tt^0$ is coarser than every triangulation $\Ttn$ ($n\in\NN$), \ie the set of edges and nodes of $\Tt^0$ is included in those of $\Ttn$. Moreover we assume that the problem~$(\LA_0)$ is feasible. Then we have of course $\min_{\HB} F_n\leq \min_{\HB} F_0<+\infty$.

Our goal is to prove that the minimizers of~\eqref{eq-pbm-discret} are close to minimizers of~\eqref{eq-prelim-problem} as the triangulation becomes thin ($n\to +\infty$ so that $h_n\to 0$). To this aim we prove the $\Gamma$-convergence of those problems towards~\eqref{eq-prelim-problem}. We refer the reader to~\cite{Dalmaso1993,braides2002gamma} for further details on $\Gamma$-convergence.

\begin{remark}
  The space $\HB$ endowed with the weak topology is a topological vector space which does not satisfy the first axiom of countability (\ie the existence of a countable base of neighborhoods at each point). Therefore, to deal with $\Gamma$-convergence with the $\HB$ weak topology, one should \textit{a priori} use the general definition of $\Gamma$-convergence (see~\cite[Def.~4.1]{Dalmaso1993}) that is valid in any topological space. 
 
   However the bounded subsets of $\Mm_b(\Omega)$ endowed with the weak-* topology are metrizable, and so are the bounded sets of $\HB$ endowed with the weak topology.
Since we are interested in the minimizers of $F_n$ ($n\in\NN$) and $F$, we might as well focus on a particular subset of $\HB$, for instance
\begin{align}
  X\eqdef \enscond{u\in \HB}{F(u)\leq \left(\min_{\HB} F_0\right)+ 1}
\end{align}
which is a bounded subset of $\HB$ by Lemma~\ref{lem-poincare}. Observe that $X$ is also compact for the weak topology of $\HB$.

Since $X$ is metrizable we may now use the following definition of $\Gamma$-convergence which is quite convenient since it is formulated in sequential terms.
\end{remark}

\begin{defn}
  We say that the sequence $(F_n)_{n\in\NN}$ $\Gamma$-converges towards $F$ for the $\HB$ weak topology if the following two properties hold:
\begin{description}
  \item[(Liminf inequality)] 
For any $u\in X$ and any sequence $(u_n)_{n\in \NN}$ such that $u_n\in X\cap \Im \ITn$ and $u_n$ converges to $u$ for the $\HB$ weak topology, 
\begin{align}
F(u)\leq \liminf_{n\to +\infty} F_n(u_n).
\label{eq-gliminf}
\end{align}

  \item[(Limsup inequality)] 
    For all $u\in X$, there exists a sequence $(u_n)_{n\in\NN}\in X\cap \Im \ITn$
 which converges towards $u$ for the $\HB$ weak topology, such that
\begin{align}
F(u)\geq \limsup_{n\to +\infty} F_n(u_n).
\label{eq-glimsup}
\end{align}
\end{description}
\label{def-gamma-cv}
\end{defn}

We are now in position to state the main result of this paper.

\begin{theorem}[$\Gamma$-convergence]
  The sequence $(F_n)_{n\in\NN}$ $\Gamma$-converge towards $F$ for the $\HB$ weak topology, and
\begin{align}
  \min_{\HB} F = \lim_{n\to +\infty} \min_{\HB} F_n.\label{eq-gamma-cv}
 \end{align}
  Moreover, every sequence $(u_n)_{n\in \NN}$ of minimizers of $F_n$ admits cluster points (for the weak topology). Each cluster point of $(u_n)_{n\in\NN}$ is a minimizer of~$F$.
  \label{thm-gamma-cv}
\end{theorem}

\begin{proof}
  Since $F_n$ coincides with $F$ on $X\cap \Im \ITn$, the liminf inequality is a straightforward consequence of the lower semi-continuity of $J$ and the continuity of $\gamma_0$ and $L$ with respect to the weak topology of $\HB$.

  Let us focus on the limsup inequality. First, we prove the result for $v\in C^3(\overline{\Omega})$. By Proposition~\ref{prop-cv-interpolation}, we know that $\ITn v$ converges towards $v$ for the strong topology of $\HB$ (since $h_n\to 0$). Hence, $\lim_{n\to +\infty} \langle L,\ITn v\rangle=1$ and the sequence $v_n\eqdef  \frac{1}{\langle L, \ITn v\rangle} \ITn v$ also converges towards $v$. By continuity of $J$ (for the intermediate topology) we have therefore $\lim_{n\to+\infty} J(v_n)=J(v)$, so that $\lim_{n\to +\infty} F_n(v_n)=F(v)$.
  For a general $u\in \HB$, we use Proposition~\ref{prop-approx-hb} to find a sequence of functions $(v_k)_{k\in\NN}\in (C^3(\overline{\Omega}))^\NN$ such that $|J(u)-J(v_k)|\leq 2^{-k}$ and $v_k$ converges towards $u$ for the intermediate topology. Applying the result to each $v_k$, we obtain a family $(v_{k,n})_{(k,n)\in\NN^2}$ and we conclude by a diagonal argument to obtain the limsup inequality.

  Thus, the functionals $F_n$ $\Gamma$-converge towards $F$. Since $X$ is compact, every sequence of minimizers has cluster points in $X$. The sequence $(F_n)_{n\in\NN}$ is equicoercive from the inequality $F\leq F_n$, and from~\cite[Theorem~7.8 and Corollary~7.20]{Dalmaso1993}, we obtain inequality~\eqref{eq-gamma-cv} and the fact that each cluster point of $(u_n)_{n\in\NN}$ is a minimizer of $F$.
\end{proof}

\begin{remark}
As the proof shows, the result not only holds for the Lagrange and Hermite elements mentioned above, but also for any affine equivalent finite element family which involves derivatives up to the second order.
\end{remark}


\section{Numerical Illustration}
\label{sec-numerics}

This section is devoted to some numerical illustrative examples. The continuous problem is discretized using either $\Pp_2$ Lagrange or $\Pp_3$ Hermite triangular finite elements and the discrete minimization problem is formulated following the method described in \cite{bleyer2013non-conforming} as a second-order conic program and solved using the dedicated software package \textsc{Mosek} \cite{mosek}. It is worth noting that the results for the $\Pp_2$ Lagrange element have been improved compared to those presented in \cite{bleyer2013non-conforming} due to the fixing of an error present in the initial numerical code. Hence, contrary to what was observed in this earlier work, the theoretical convergence result obtained in the present paper is observed numerically.

\begin{figure}
\begin{center}
\includegraphics[width = 0.8\textwidth]{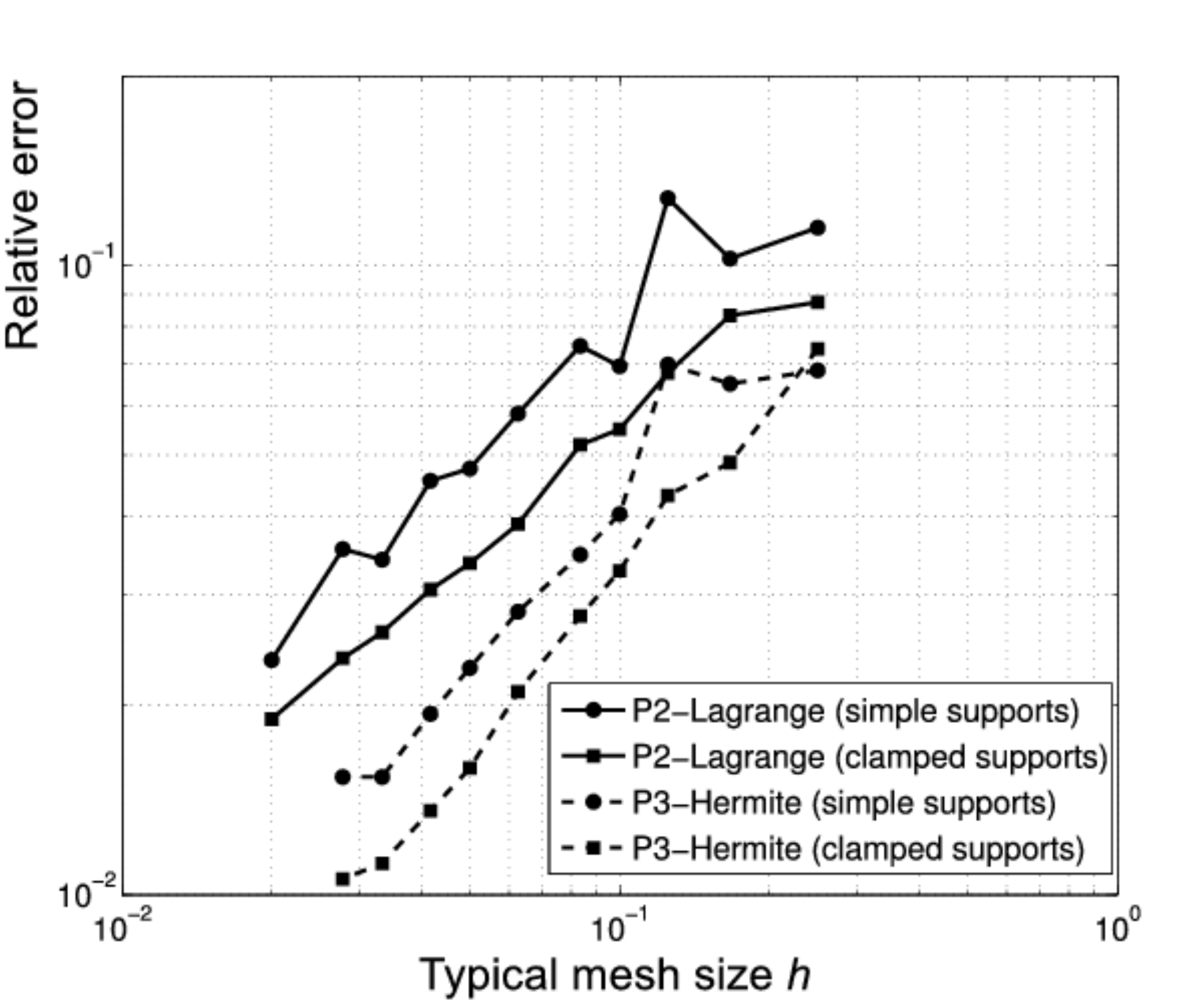} 
\end{center}
\caption{Relative error versus mesh size for the square plate problem for the $\Pp_2$ Lagrange and $\Pp_3$ Hermite element for different boundary conditions}
\label{rel-error-square}
\end{figure}

\begin{figure}
\begin{center}
\includegraphics[width = \textwidth]{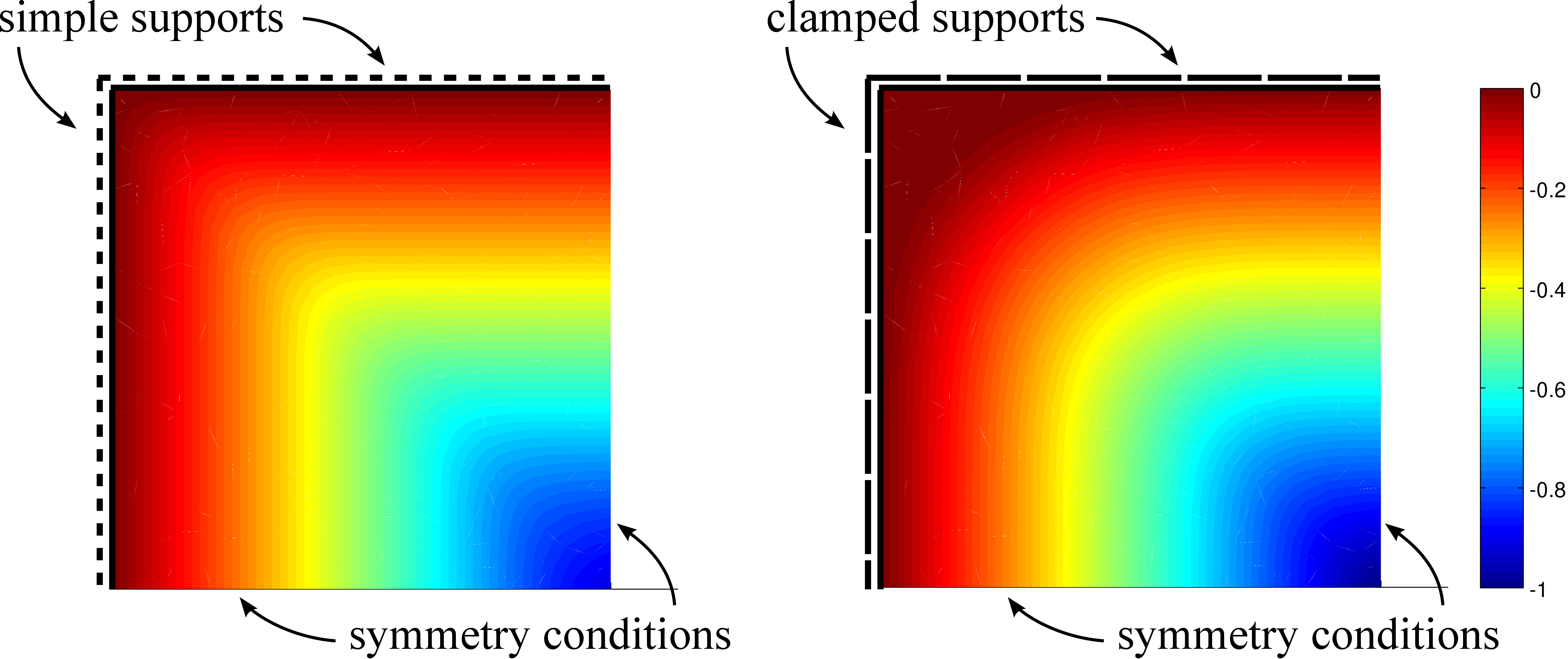} 
\end{center}
\caption{Representation of the obtained failure mechanism $u$ for the square plate problem with simple supports (left) and clamped supports (right) using the $\Pp_2$ Lagrange element\\}
\label{fail-mech}
\end{figure}

\subsection{Homogeneous Case}

The first example considers the problem of a square plate of side $a$ under a uniform transversal reference loading $L(x)=L$. The plate is supposed to obey the Johansen strength criterion with $M_0^+=M_0^-=M_0$. The boundary conditions are either simple supports (i.e. $\Gamma_D = \partial \Omega$, $\Gamma_N = \emptyset$) or clamped supports (i.e. $\Gamma_D = \Gamma_N = \partial \Omega$). In each case, an analytical solution is available for the ultimate load : $\lambda^+ = 24 M_0/(La^2)$ for the simple supports \cite{save1997plastic} and $\lambda^+ = 42.851 M_0/(La^2)$ for the clamped supports \cite{fox1974limit}. It is worth noting that the solution for simple supports is very simple as the optimal velocity field corresponds to four parts separated by the plate diagonals which rotate along the four boundaries. For the clamped problem, the optimal velocity field is much more complicated as it consists of a combination of developable and anticlastic surfaces, together with undeflected corners \cite{fox1974limit}.

Only a quarter of the plate has been discretized and proper symmetry conditions ($u$ free and $\frac{\partial u}{\partial \nu} = 0$) have been imposed on the axes of symmetry. Different upper bound estimates $\lambda_h$ of the exact ultimate load have been obtained when varying the typical size $h$ of a finite element which ranged between $0.25a$ and $0.02a$. The relative errors $\frac{\lambda_h-\lambda^+}{\lambda^+}$ have been represented in figure \ref{rel-error-square} for both types of boundary conditions and finite elements. It can be observed that the discrete estimates seem to converge to the exact value with approximately the same rate for both boundary conditions but the $\Pp_3$ Hermite element exhibits a higher convergence rate than the $\Pp_2$ Lagrange element.

It is also worth noting that the convergence is much smoother for the clamped problem than the simple supports problem. This can be attributed to the fact that the optimal failure mechanism consists of a discontinuity of the gradient along the diagonal, the quality of the solution is, therefore, much more mesh dependent since we are using an unstructured mesh, the edges of which are not aligned along the diagonal \textsl{a priori}. Finally, optimal discrete velocity fields have been plotted for both problems in figure \ref{fail-mech}. It can be observed that the mechanism for simple supports seems to reproduce a mechanism with a concentrated rotation discontinuity along the diagonal and rigid parts rotating along the plate boundary. For the clamped case, the mechanism is indeed more complicated and an undeflected region is observed near the plate corner.

\begin{figure}
\begin{center}
\includegraphics[width = 0.8\textwidth]{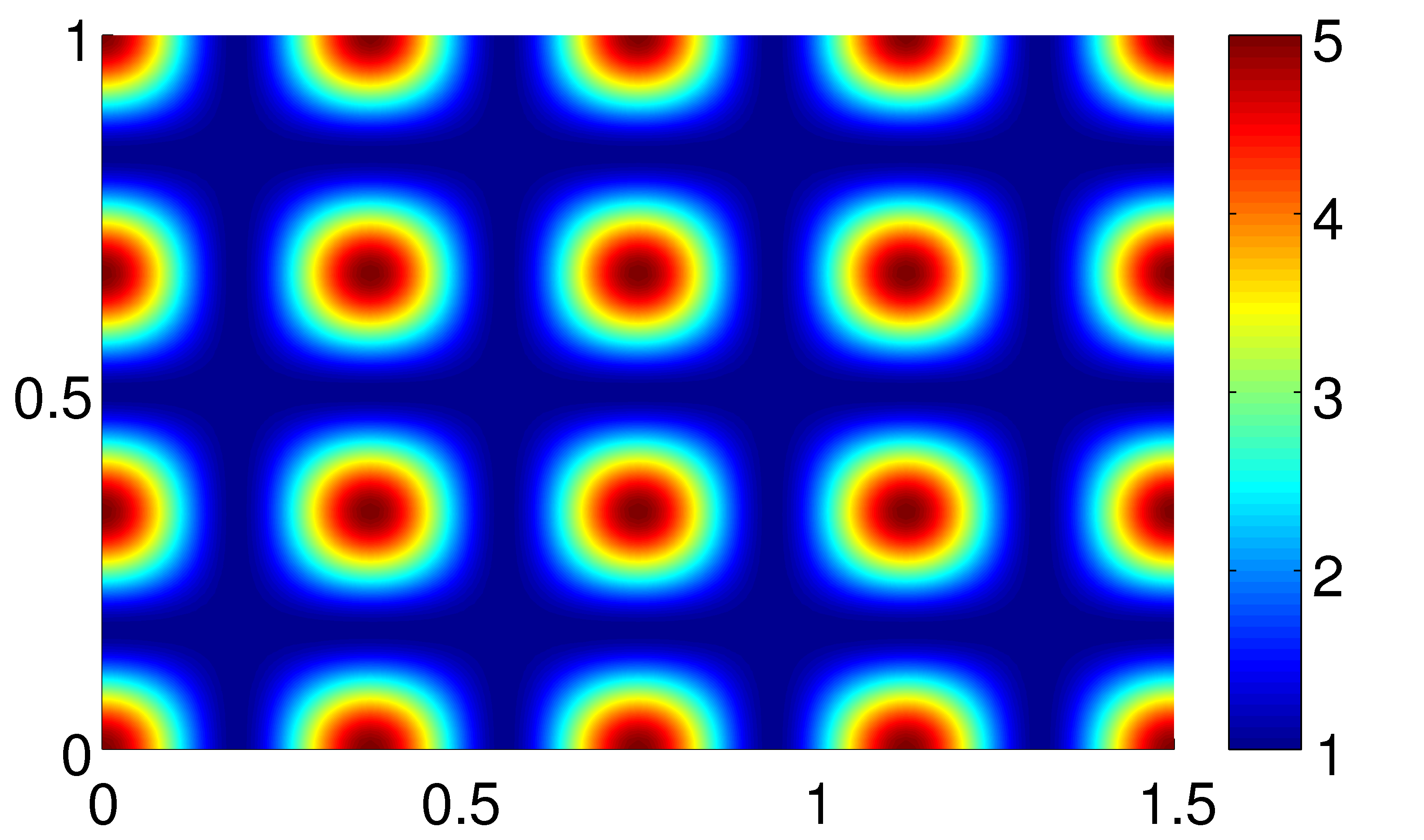} 
\end{center}
\caption{Representation of the variation of the uniaxial strength in bending $M_0(x_1,x_2)$ for the non-uniform example\\}
\label{strength-variation}
\end{figure}

\subsection{Inhomogeneous Case}

The second example involves a rectangular plate domain $\Omega=[0;1.5]\times [0;1]$ with simple supports along the plate boundary and subjected to a uniform transversal loading $L$. In this example, the case of a non-uniform distribution of the penalty $\pi$ function is considered. More particularly, it has been assumed that the plate obeys a von Mises strength criterion with a strength in uniaxial bending $M_0$ which is non-uniform throughout the plate domain $\Omega$ and which has been taken as $M_0(x_1,x_2) = (\cos(\frac{16\pi}{3}x_1)+1)(\cos(6\pi x_2)+1)+1$ (see Figure \ref{strength-variation}). In this case, the penalty function can be written for all $x\in \Omega$ as $\pi(x,D^2u) = M_0(x_1,x_2)\Pi(D^2u)$ where $\Pi$ (which does not depend explicitly on the point $x$) is the support function corresponding to the von Mises criterion with a strength in bending parameter of unit value.

\begin{figure}
\begin{center}
\subfloat[Uniform distribution of $M_0$]{\label{mech-uniform}
\includegraphics[width = 0.8\textwidth]{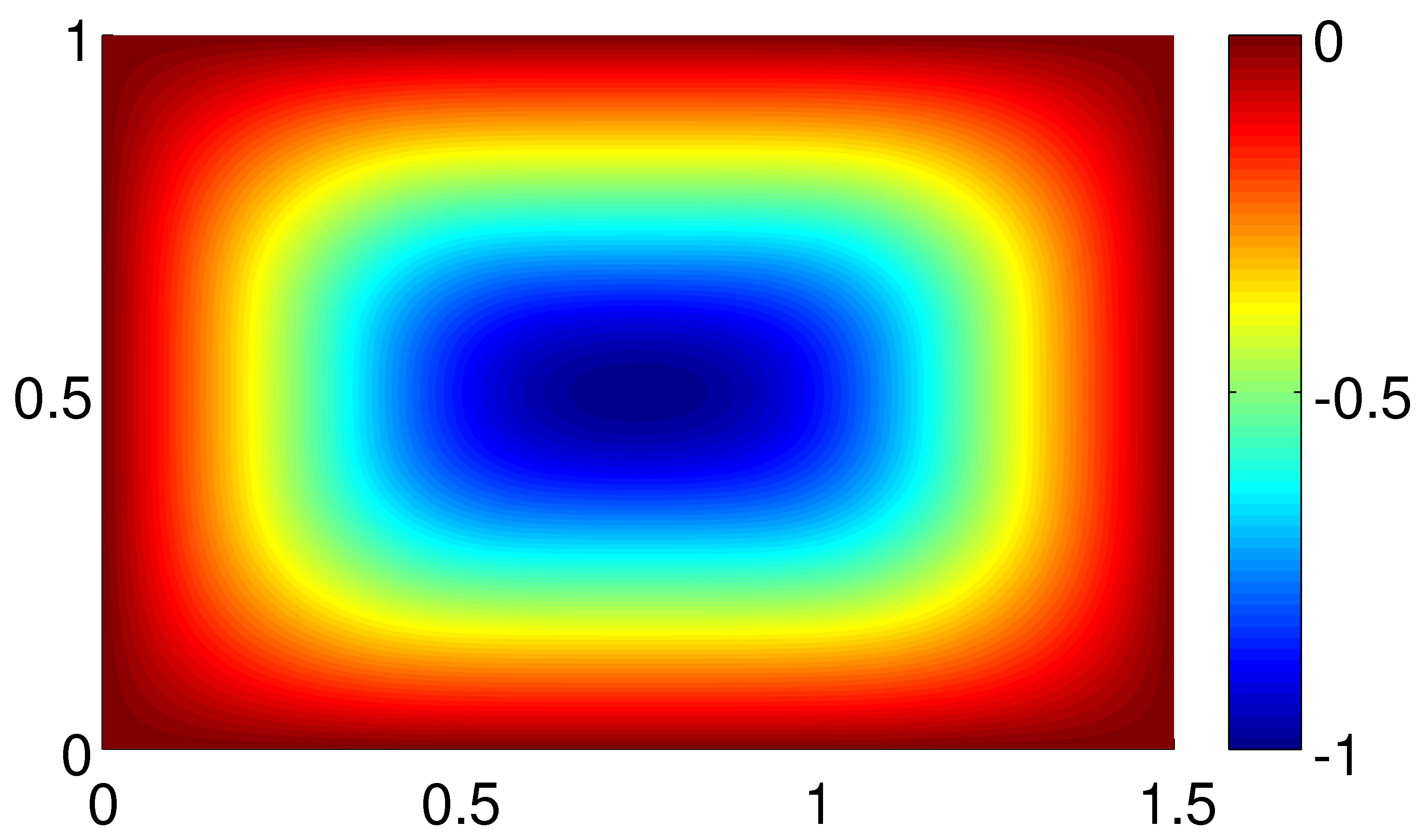} 
}
\\
\subfloat[Non-uniform distribution of $M_0$]{\label{mech-nonuniform}
\includegraphics[width = 0.8\textwidth]{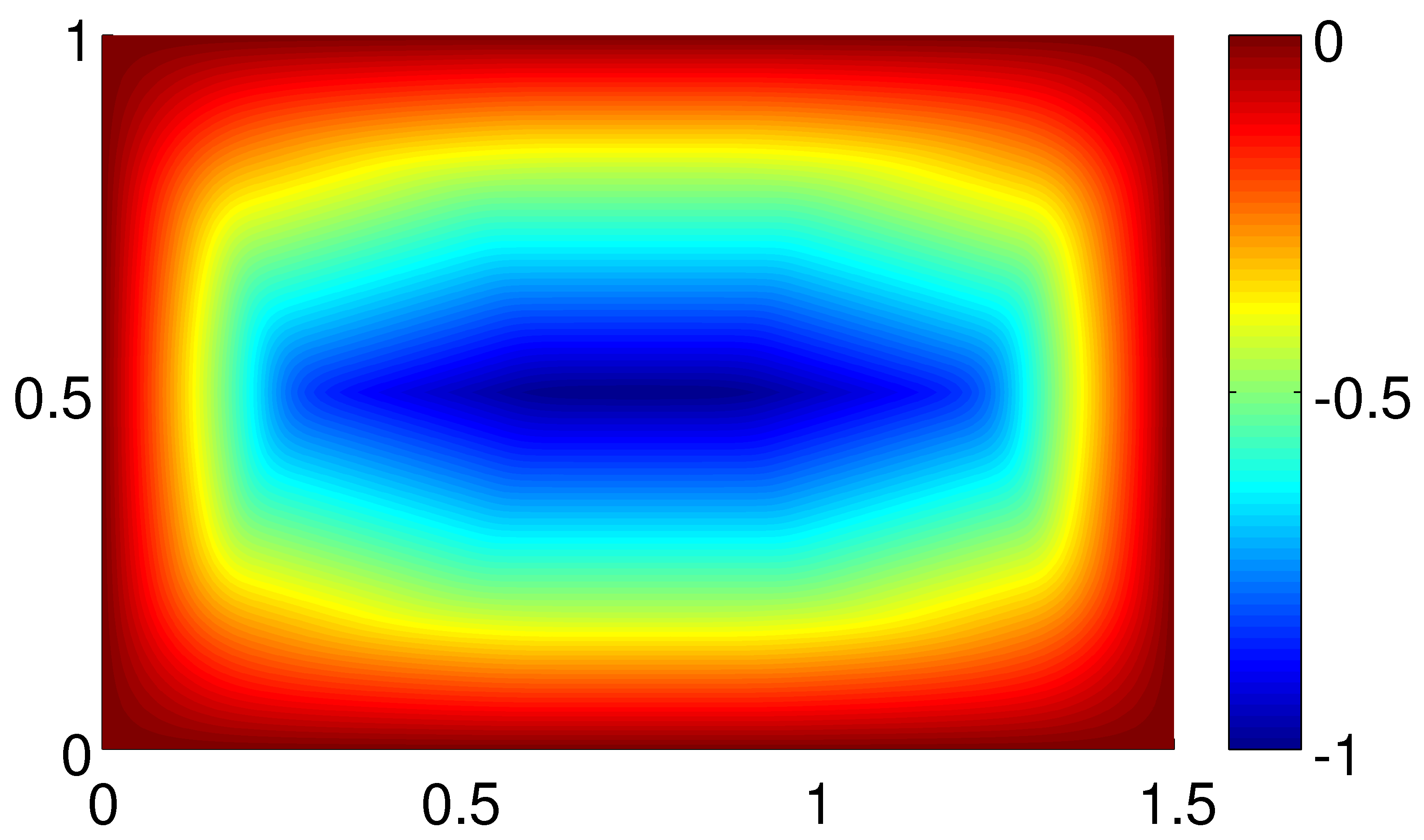} 
}
\end{center}
\caption{Comparison of the relative distribution of the optimal failure mechanism obtained for a non-uniform and a uniform distribution of the strength parameter\\}
\label{strength-variation}
\end{figure}

The optimal failure mechanism obtained for this problem has been represented in figure \ref{mech-nonuniform} whereas the optimal mechanism obtained for the same problem with a uniform distribution of the strength parameter $M_0=1$ has been represented for comparison in figure \ref{mech-uniform}. One can observe that both mechanisms exhibit some distinctive features since the contour lines obtained in the uniform case seem to be quite smooth whereas they seem to be piecewise linear for the non-uniform case, at least in the region close to the center of the plate.

\begin{figure}
\begin{center}
\subfloat[Uniform distribution of $M_0$]{\label{pi-uniform}
\includegraphics[width = 0.8\textwidth]{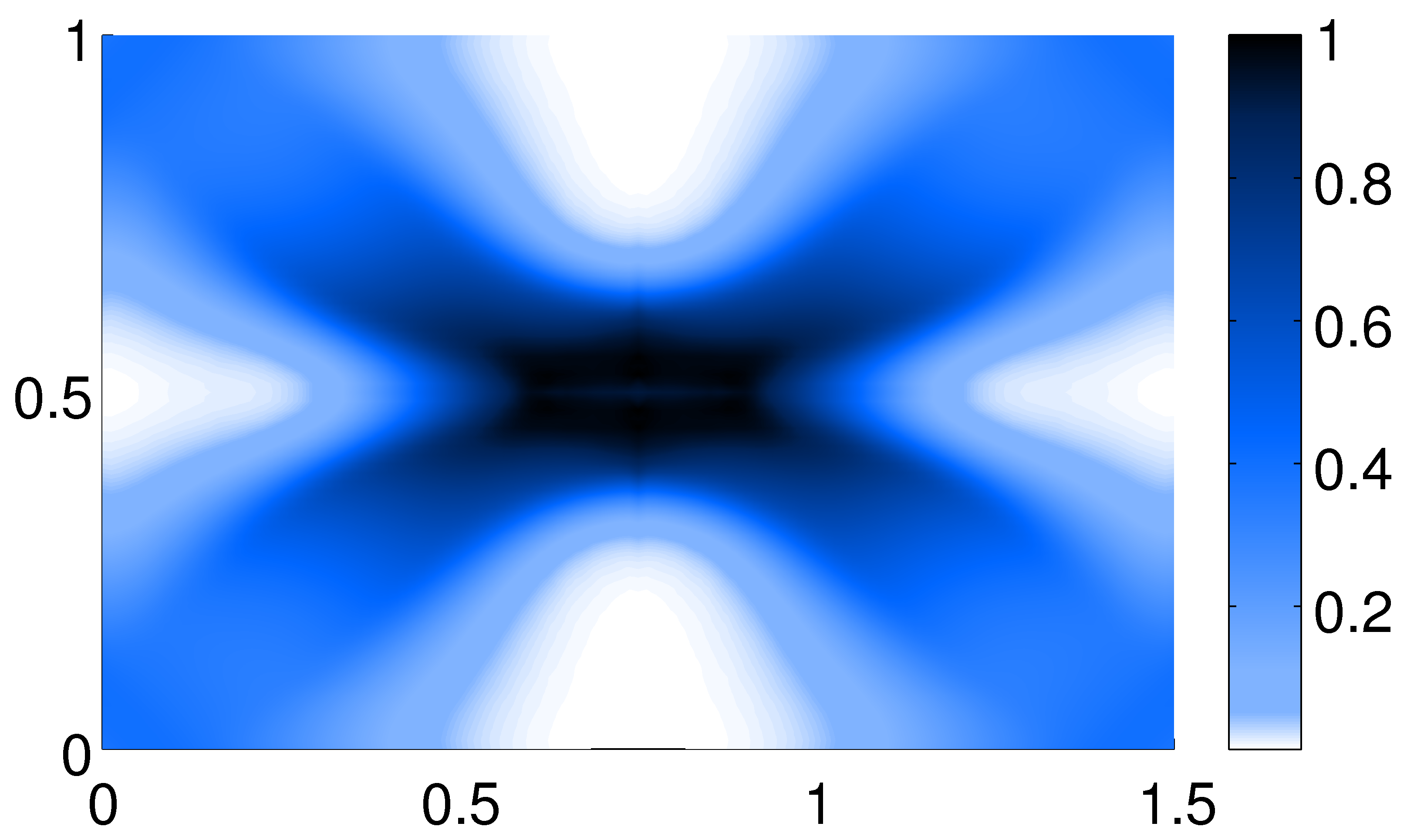} 
}
\\
\subfloat[Non-uniform distribution of $M_0$]{\label{pi-nonuniform}
\includegraphics[width = 0.8\textwidth]{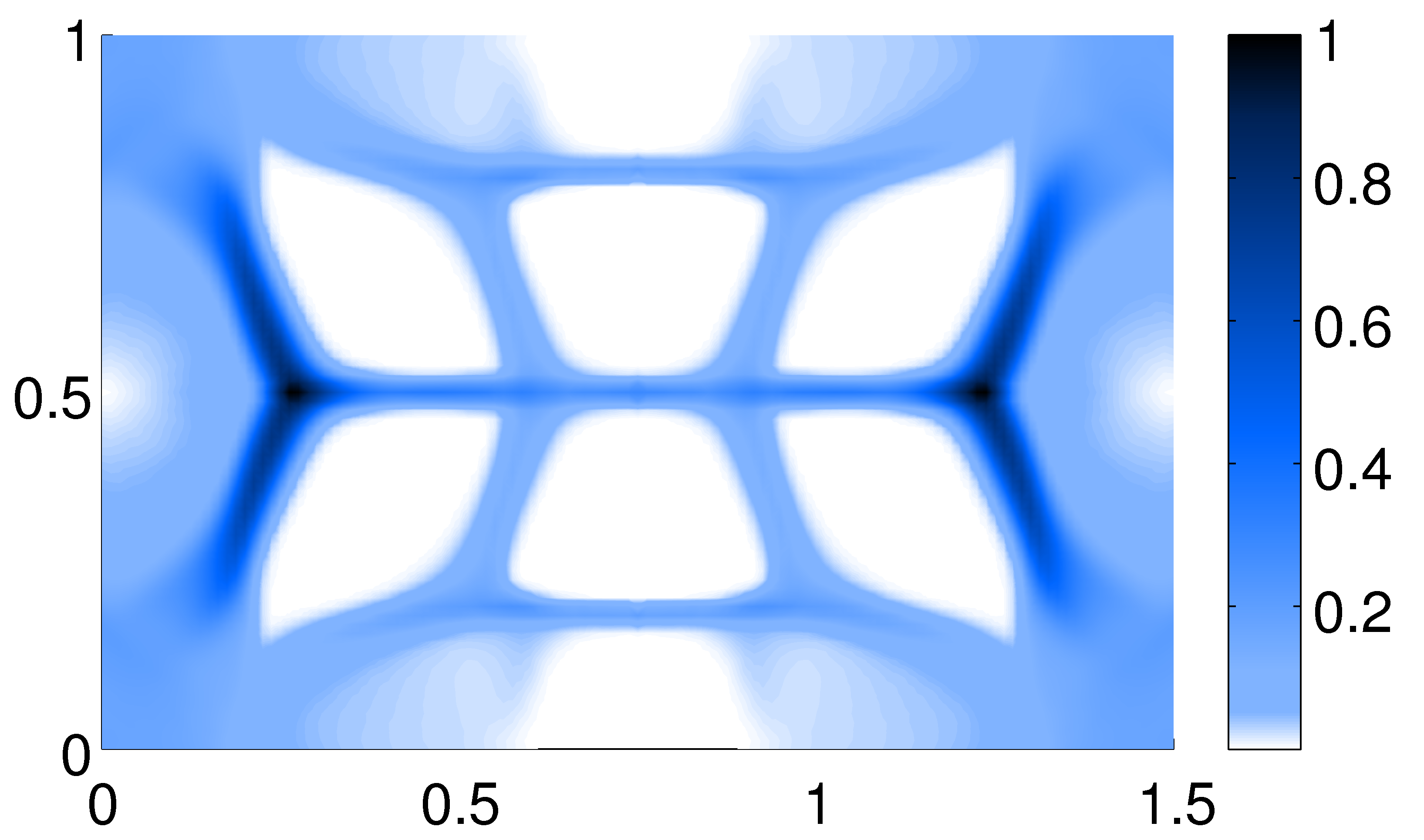} 
}
\end{center}
\caption{Comparison of the relative distribution of $\Pi(D^2u(x))$ obtained for a non-uniform and a uniform distribution of the strength parameter\\}
\label{strength-variation}
\end{figure}

This observation can be further interpreted by looking at the relative distribution of the function $\Pi(D^2u(x))$ for the optimal mechanism in both cases. In the uniform case (Figure \ref{pi-uniform}), one can observe that the quantity $\Pi(D^2u(x))$, which is some kind of Euclidean norm of the curvature of the optimal mechanism, is non-zero in broad regions around the plate diagonals. Besides, there is no narrow region of highly concentrated curvature meaning that the rotation field (the velocity gradient) does not exhibit any sharp features. On the contrary, in the non-uniform case (Figure \ref{pi-nonuniform}), it can be seen that there exists zones of highly localized curvature deformation situated along the minima of the $M_0$ distribution. In the regions where $M_0$ attains its maximum, there is on the contrary no deformation at all (white regions). Therefore, in these regions, the optimal field is linear and there is a discontinuity of the rotation field where $M_0$ is minimal. This is valid essentially at the center of the plate since it is not possible to obtain a piecewise linear velocity field along the minima around the corners due to the boundary conditions.

\section*{Conclusion}

In this article, we have provided for the first time a rigorous analysis of the convergence of second order finite element discretizations for the limit analysis of thin plates. This requires a careful definition of the corresponding continuous problem over the space of Bounded Hessian functions. Handling second order derivatives which are measures requires some special approximation arguments to ensure the convergence of the finite element method. 

Let us emphasize that, although we have insisted on the Lagrange and Hermite finite elements, the proposed result holds in the more general case of any finite element system which is affine equivalent (see~\cite{brenner2008}) to some reference element an which involves derivatives up to the second order. 
Future works may extend this approach to finer models such as the study of thick plates in bending and shear forces.

\section*{Acknowledgements} 

The work of Gabriel Peyr\'e and Vincent Duval has been supported by the European Research Council (ERC project SIGMA-Vision).

\bibliographystyle{alpha}
\bibliography{biblio}

\end{document}